\documentclass[10pt]{ijnam}
\theoremstyle{plain}
\newtheorem{mytheo}{Theorem}[section]
\newtheorem{mydef}[mytheo]{Definition}

\newtheorem{assum}[mytheo]{Assumption}

\newtheorem{rem}[mytheo]{Remark}

\usepackage{graphicx}%
\usepackage{multirow}%
\usepackage{amsmath,amssymb,amsfonts}%
\usepackage{amsthm}%
\usepackage{mathrsfs}%
\usepackage{xcolor}%
\usepackage{algpseudocode}%
\usepackage{listings}%
\usepackage{color}%
\usepackage{subfigure}
\copyrightinfo{2004}{} 

\usepackage{mydef}

\begin{document}

\title[]
{Two new families of fourth-order explicit exponential Runge--Kutta methods with four stages for first-order
differential systems}

\author[X. Hu]{Xianfa Hu}
\address{
Department of Mathematics, Shanghai Normal University, 
 Shanghai, 200234, P.R.China
}
\email{zzxyhxf@163.com}

\author[Y. Fang]{Yonglei Fang}
\address{School of Mathematics and Statistics, Zaozhuang University,
 Zaozhuang, 277160,  P.R.China}
 \email{ylfangmath@163.com}
 
 \author[B. Wang]{Bin Wang}
 \address{School of Mathematics and Statistics, Xi'an Jiaotong University, 
 Xi'an, 710049, P.R.China}
  \email{wangbinmaths@xjtu.edu.cn}

\subjclass[2000]{65L05, 65L06}

\abstract{In this paper,  two new families of fourth-order explicit
exponential Runge--Kutta (ERK)  methods  with four stages    are
studied for solving  first-order differential systems $y'(t)+My(t)=f(y(t))$.  By comparing the  Taylor series of  the exact solution, the order conditions
 of these ERK methods are derived, which are exactly identical to the order
 conditions of  explicit  Runge--Kutta methods, and  these ERK methods reduce to  classical
 Runge--Kutta methods once $M\rightarrow \mathbf{0}$.
Moreover,  we analyze the stability properties and  the convergence of the new  methods. Several numerical examples are implemented  to illustrate  the  accuracy and efficiency of these ERK methods by comparison with standard exponential integrators.
}

\keywords{Exponential Runge--Kutta methods, first-order differential equations, the order conditions, and  convergence}

\maketitle
\section{Introduction}

Classical Runge--Kutta (RK) methods are extensively
recognized by researchers and engineers for its simple idea
and concise expression \cite{Butcher2008,Qi2012,Wang2021,Zhang2014}, and standard ERK methods with the stiff-order conditions (comprise the classical order conditions) were formulated by Hochbruck et al. \cite{Hochbruck2005a,Hochbruck2005b,Hochbruck2010}, which can be viewed as  an extension of classical RK methods.
The idea of exponential integrators is primarily concerned with the use of the variation-of-constants formula (see \cite{Hochbruck2010}), normally,   the coefficients of  exponential integrators
 are  exponential functions of  the underlying matrix in the literature, this fact means that 
  the implementation of  exponential integrators requires the
evaluations of matrix exponentials and vectors. In order to reduce the computational cost,
 two new kinds of ERK methods up to order three were formulated in \cite{Wang2022}.
As a sequel to this work, we study two new
 families of 
  fourth-order explicit ERK methods with four stages for the
first-order differential system in this paper, which are termed the modified and  simplified versions
  of  fourth-order explicit ERK  integrators, respectively.

 In this paper, we consider the evolution equation  and ordinary  differential equation of the form
\begin{equation}\label{equ1}\left\{
\begin{array}{l}
y'(t)+My(t)=f(y(t)),\quad t\in[t_0,t_{\rm
end}],\cr\noalign{\vskip1truemm} y(t_0)=y_0,
\end{array}
\right.\end{equation} {where} the matrix  $M\in \mathbb R^{m\times
m}$ is a symmetric positive definite or skew-Hermitian. Problems of the form (\ref{equ1})
arise frequently in a variety of applied science  such as quantum mechanics,
flexible mechanics, and semilinear parabolic problems. 
  
 Exponential integrators have received much attention due to its high accuracy and stability (see, e.g., \cite{Cell2008,Hochbruck1998,Lambert1972,Lawson1967,Li2016,Mei2022,Nosett1969,Pope1963,Wu2019,Wu2018,Wu2021}), and  the product of a matrix exponential and a vector  is effectively evaluated by the Krylov subspace method \cite{Eiermann2006,Hochbruck1997} or a Pade approximant (see, e.g., \cite{Berland2007,Moler2003}). In recent years, Fang et al. \cite{Fang2021} presented explicit pseudo two-step exponential Runge-Kutta methods for first-order differential equations, Du et al. \cite{Du2019,Du2021} considered the exponential time differencing method  for semilinear parabolic problems, and Li et al. \cite{Li2020} formulated exponential cut-off methods for preserving maximum principle.  It should be noted that ERK methods are based on the stiff-order conditions.   However, as stated by Berland et al. in \cite{Berland2005}, the stiff-order conditions are relatively strict. It is true that the fourth-order explicit ERK method \eqref{Tableau-ERK51}
with five stages was proposed in \cite{Hochbruck2005a}, which only satisfied the stiff-order conditions
in a weak form. Here,  the coefficients of our fourth-order explicit ERK methods are real constants and the
study is related to the classical order conditions.

 The outline of the paper is as follows. In Section \ref{sec2}, we  derive the order conditions of the modified version of fourth-order explicit ERK  methods. Section \ref{sec3} is devoted to presenting  the simplified version of  explicit ERK  methods of order four. The linear stability properties and  the convergence of our explicit ERK methods  are analyzed in Section \ref{sec4}. Numerical results  present the accuracy and efficiency of these  ERK methods when  applied to the averaged system in wind-induced oscillation,  the H\'{e}non-Heiles Model, the Allen-Cahn equation, the sine-Gorden equation and the nonlinear Schr\"{o}dinger equation in Section \ref{sec5}. The concluding remarks are included in the last section.

\section{ A modified version of fourth-order explicit ERK methods}\label{sec2}

The internal stages  and the  update of the modified version of ERK methods inherit and modify the form of classical RK methods, respectively.
\begin{mydef}\label{definition1} (\cite{Wang2022}) An $s$-stage modified version of  exponential
{Runge--Kutta} (MVERK) method for the numerical integration
(\ref{equ1}) is defined as
\begin{equation}\label{MVERK}
\left\{
\begin{array}{l}
\displaystyle Y_i=y_0+h\sum\limits_{j=1}^sa_{ij}(-MY_j+f(Y_j)), \quad
i=1,\ldots,s,\cr\noalign{\vskip1truemm}
\displaystyle y_{1}=e^{-hM}y_0+h\sum\limits_{i=1}^sb_{i}f(Y_i)+w_s(z),
\end{array}
\right.
\end{equation}
{where} $a_{ij} $, $b_i$
 are real constants for {$i,j=1,\ldots, s$},
 $Y_i\approx y(t_0+c_ih)$ for $i=1,\ldots,s$,  $w_s(z)$ is a suitable matrix-valued function of  $hM$,  and $w_s(z)\rightarrow 0$ when $M\rightarrow \mathbf{0}$. \end{mydef}

 In  \cite{Wang2022}, it is  clearly indicated that $w_s(z)$ is independent of the  matrix exponential.  In fact,  the $w_s(z)$ is also related to the term $f(\cdot)$ and initial value $y_0$, and  MVERK methods with the same order share the same $w_s(z)$. Especially, if we consider the MVERK method with order one, then $w_s(z)=0$.  From \eqref{MVERK}, it is  clear that the MVERK method can exactly integrate the first-order homogeneous linear system $y'(t)=-My(t), \ y(0)=y_0,$
which has the exact solution
$y(t)=e^{-tM}y_0.$ The property of the method (\ref{MVERK}) is significant. For linear  oscillatory problems, the exponential contains the full information on linear oscillations in contrast to classical numerical methods (non-exponential).  The method (\ref{MVERK}) can be displayed by the following Butcher Tableau
\begin{equation}\label{Tableau-MVERK4}
\begin{aligned} &\quad\quad\begin{tabular}{c|c|c}
 ${c}$&$\mathbf{I}$&${A}$ \\
 \hline
  $\raisebox{-1.3ex}[1.0pt]{$e^{-hM}$}$ & $\raisebox{-1.3ex}[1.0pt]{$w_s(z)$}$&$\raisebox{-1.3ex}[1.0pt]{${b}^{\mathrm{T}}$}$  \\
\end{tabular}
~=
\begin{tabular}{c|c|ccc}
 ${c}_1$&$I$&${a}_{11}$&$\cdots$&${a}_{1s}$\\
$\vdots$& $\vdots$ & $\vdots$&$\vdots$&$\vdots$\\
 ${c}_s$ &$I$&  ${a}_{s1}$& $\cdots$& ${a}_{ss}$\\
 \hline
 $\raisebox{-1.3ex}[1.0pt]{$e^{-hM}$}$&$\raisebox{-1.3ex}[1.0pt]{$w_s(z)$}$&$\raisebox{-1.3ex}[1.0pt]{${b}_1$}$&\raisebox{-1.3ex}[1.0pt]{$\cdots$} &  $\raisebox{-1.3ex}[1.0pt]{${b}_s$}$\\
\end{tabular}
\end{aligned}
\end{equation}
with  $c_i=\sum\limits_{j=1}^sa_{ij}$  for $i=1,\ldots,s$.

In view of  \eqref{MVERK}, the internal stages of the MVERK method are independent of matrix exponentials, and the update remains some properties of the matrix exponential.  Once $M \rightarrow \mathbf{0}$ i.e., $y^{\prime}(t)=f(y(t))$, then $e^{-hM}=I$, $ w_s(z)=\mathbf{0}$, and the MVERK method reduces to a classical RK method
\begin{equation}\label{RK}
\left\{
\begin{array}{l}
Y_i=y_0+h\sum\limits_{j=1}^sa_{ij}f(Y_j), \
i=1,\ldots,s,\cr\noalign{\vskip3truemm}
y_{1}=y_0+h\sum\limits_{i=1}^sb_{i}f(Y_i).
\end{array}
\right.
\end{equation}
Therefore,  MVERK methods  can be considered as an extension of classical RK methods.

  A numerical method is said to be of order $p$ if the Taylor expansion of a numerical solution $y_1$ and the exact solution $y(t_0+h)$ coincide up to $h^p$ about $y_0$. Under the assumption $y(t_0)=y_0$,  we denote $g(t_0)=-My(t_0)+f(y(t_0))$, and the Taylor series of $y(t_0+h)$ is given by
\small{\begin{align*}
\displaystyle &y(t_0+h)=y(t_0)+hy'(t_0)+\frac{h^2}{2!}y''(t_0)+\frac{h^3}{3!}y'''(t_0)+\frac{h^4}{4!}y^{(4)}(t_0)+\mathcal{O}{(h^5)}\cr\noalign{\vskip3truemm} &=y(t_0)+hg(t_0)+\frac{h^2}{2!}(-M+f'_y(y_0))g(t_0)+\frac{h^3}{3!}\Big(M^2g(t_0)+(-M+f'_y(y(t_0)))f'_y(y(t_0))g(t_0)
\cr\noalign{\vskip2truemm}
& \quad -f'_y(y(t_0))Mg(t_0) +f''_{yy}(y(t_0))(g(t_0),g(t_0))\Big)+\frac{h^4}{4!}\Big(-M^3g(t_0)+M^2f'_y(y(t_0))g(t_0)\cr\noalign{\vskip3truemm}
&\quad -Mf'_y(y(t_0))(-M+f'_y(y(t_0)))g(t_0)
-Mf''_{yy}(y(t_0))(g(t_0),g(t_0))+f'_y(y(t_0)) (-M\cr\noalign{\vskip2truemm}
&\quad+f'_y(y(t_0))) (-M+f'_y(y(t_0)))g(t_0)+f'_y(y(t_0))f''_{yy}(y(t_0))(g(t_0),g(t_0)) +f'''_{yyy}(y(t_0))(g(t_0),\cr\noalign{\vskip3truemm}
&\quad g(t_0),g(t_0))+3f''_{yy}(y(t_0))\big((-M+f'_y(y(t_0)))g(t_0),g(t_0)\big)\Big)+\mathcal{O}{(h^5)}.
\end{align*}}
Before we derive the order conditions of  fourth-order explicit MVERK methods with four stages, it is necessary to present the classical order conditions  of RK methods:
\begin{equation}\label{order condition}
\begin{aligned}
 \includegraphics{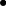} \qquad &b_1+b_2+b_3+b_4&=1, \cr\noalign{\vskip2truemm}
\includegraphics{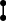} \qquad &b_2c_2+b_3c_3+b_4c_4&=\frac{1}{2},\cr\noalign{\vskip2truemm}
\includegraphics{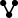} \qquad  &b_2c_2^2+b_3c_3^2+b_4c_4^2&=\frac{1}{3},\cr\noalign{\vskip2truemm}
\includegraphics{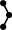} \qquad &b_3a_{32}c_2+b_4a_{42}c_2+b_4a_{43}c_3&=\frac{1}{6},\cr\noalign{\vskip2truemm}
 \includegraphics{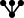} \qquad&b_2c_2^3+b_3c_3^3+b_4c_4^3&=\frac{1}{4},\cr\noalign{\vskip2truemm}
 \includegraphics{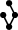} \qquad&b_3c_3a_{32}c_2+b_4c_4a_{42}c_2+b_4c_4a_{43}c_3&=\frac{1}{8},\cr\noalign{\vskip2truemm}
 \includegraphics{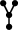} \qquad&b_3a_{32}c_2^2+b_4a_{42}c_2^2+b_4a_{43}c_3^2&=\frac{1}{12},\cr\noalign{\vskip2truemm}
 \includegraphics{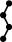} \qquad&b_4a_{43}a_{32}c_2&=\frac{1}{24},
 \end{aligned}
\end{equation}
 with the appropriate assumption $c_i=\sum\limits_{j=1}^{i-1}a_{ij}$, for $i=1,2,3,4$.  It is  very clear that the equation  (\ref{order condition}) has infinitely many solutions.

 The following theorem will present the order conditions
 of fourth-order explicit MVERK methods with four stages are identical to  \eqref{order condition}.

\begin{mytheo}
 If the coefficients of the four-stage explicit MVERK method with $w_4(z)$
\begin{equation}\label{MVERK4}
\left\{
\begin{array}{l}
Y_1=y_0, \cr\noalign{\vskip2truemm}
Y_2=y_0+ha_{21}(-MY_1+f(Y_1)),\cr\noalign{\vskip2truemm}
Y_3=y_0+h\Big(a_{31}(-MY_1+f(Y_1))+a_{32}(-MY_2+f(Y_2))\Big), \cr\noalign{\vskip2truemm}
Y_4=y_0+h\Big(a_{41}(-MY_1+f(Y_1))+a_{42}(-MY_2+f(Y_2))+a_{43}(-MY_3+f(Y_3))\Big), \cr\noalign{\vskip2truemm}
y_1=e^{-hM}y_0+h\Big(b_1f(Y_1)+b_2f(Y_2)+b_3f(Y_3)+b_4f(Y_4)\Big)+w_4(z),
\end{array}
\right.
\end{equation}
where
\begin{equation*}
\begin{aligned}
w_4(z)=&-\frac{h^2}{2!}Mf(y_0)+\frac{h^3}{3!}\big(M^2f(y_0)-Mf'_y(y_0)g(y_0)\big) +\frac{h^4}{4!}\Big(-M^3f(y_0)+M^2f'_y(y_0)g(y_0)\cr\noalign{\vskip2truemm}
&-Mf''_{yy}(y_0)\big(g(y_0),g(y_0)\big)
 -Mf'_y(y_0)(-M+f'_y(y_0))g(y_0)\Big)
 \end{aligned}
 \end{equation*}
and $g(y_0)=-My_0+f(y_0)$, satisfy the order conditions \eqref{order condition}, then the explicit MVERK method has order four.
\end{mytheo}
\begin{proof}
 Under the assumption  $y(t_0)=y_0$, the Taylor series of  $y_1$ is
 \small{
\begin{equation*}
\begin{aligned}
y_1=&\big(I-hM+\frac{h^2M^2}{2!}-\frac{h^3M^3}{3!}+\frac{h^4M^4}{4!}+\mathcal{O}{(h^5)}\big)y_0+h\Big(b_1f(y_0)+b_2f(y_0)+b_2a_{21}hf'_y(y_0)g(y_0)\cr\noalign{\vskip1truemm}
&+b_2a_{21}^2\frac{h^2}{2!}f''_{yy}(y_0)(g(y_0),g(y_0))+b_2a_{21}^3\frac{h^3}{3!}f'''_{yyy}(y_0)(g(y_0),g(y_0),g(y_0))+b_3f(y_0)+b_3a_{31}h\cr\noalign{\vskip2truemm}
&\cdot f'_y(y_0)g(y_0)+b_3a_{32}hf'_y(y_0)(-MY_2+f(Y_2))+b_3\frac{h^2}{2!}f''_{yy}(y_0)\big(a_{31}g(y_0)+a_{32}(-MY_2+f(Y_2)),\cr\noalign{\vskip2truemm}
&a_{31}g(y_0)+a_{32}(-MY_2+f(Y_2))\big)+b_3\frac{h^3}{3!}f'''_{yyy}(y_0)\big(a_{31}g(y_0)+a_{32}(-MY_2+f(Y_2)),a_{31}g(y_0)\cr\noalign{\vskip2truemm}
&+a_{32}(-MY_2+f(Y_2)),a_{31}g(y_0)+a_{32}(-MY_2+f(Y_2))\big)+b_4f(y_0)+b_4hf'_y(y_0)\big(a_{41}g(y_0)\cr\noalign{\vskip2truemm}
&+a_{42}(-MY_2+f(Y_2)) +a_{43}(-MY_3+f(Y_3))\big) +b_4\frac{h^2}{2!}f''_{yy}(y_0)\big(a_{41}g(y_0)+a_{42}(-MY_2+f(Y_2))\cr\noalign{\vskip2truemm}
&+a_{43}(-MY_3+f(Y_3)),a_{41}g(y_0)+a_{42}(-MY_2+f(Y_2))+a_{43}(-MY_3+f(Y_3))\big)+w_4(z)\cr\noalign{\vskip2truemm}
&+b_4\frac{h^3}{3!}f'''_{yyy}(y_0)\big(a_{41}g(y_0)+a_{42}(-MY_2+f(Y_2))+a_{43}(-MY_3+f(Y_3)),a_{41}g(y_0)+a_{42}(-MY_2
\cr\noalign{\vskip2truemm}
&+f(Y_2))+a_{43}(-MY_3+f(Y_3)),a_{41}g(y_0)+a_{42}(-MY_2+f(Y_2))+a_{43}(-MY_3+f(Y_3))\big)\Big).
\end{aligned}
\end{equation*}}

Inserting the internal stages $Y_2,\ Y_3,\ Y_4$ into  the above formula, and combining the Taylor  series of $f(Y2),\ f(Y3),\ f(Y4)$ about $y_0$, we obtain
\small{
\begin{equation*}
\begin{aligned}
y_1=&y_0-hMy_0+h(b_1+b_2+b_3+b_4)f(y_0)-\frac{h^2Mg(y_0)}{2!}+h^2\big(b_2a_{21}+b_3(a_{31}+a_{32})+b_4(a_{41}\cr\noalign{\vskip2truemm}
&+a_{42}+a_{43})\big)f'_y(y_0)g(y_0)
+\frac{h^3}{3!}M(M-f'_y(y_0))g(y_0)+h^3\big(b_3a_{32}a_{21}+b_4a_{42}a_{21}+b_4a_{43}(a_{31}\cr\noalign{\vskip2truemm}
&+a_{32})\big)f'_y(y_0)(-M+f'_y(y_0))g(y_0)+\frac{h^3}{2!}\big(b_2a_{21}^2+b_3(a_{31}^2+2a_{31}a_{32}+a_{32}^2)+b_4(a_{41}^2+a_{42}^2\cr\noalign{\vskip2truemm}
&+a_{43}^2+2(a_{41} a_{42}+a_{41}a_{43}+a_{42}a_{43}))\big)f''_{yy}(y_0)(g(y_0)g(y_0)) +\frac{h^4}{4!}\big(-M^3g(y_0)+M^2f'_y(y_0)\cr\noalign{\vskip2truemm}
&\cdot g(y_0)-Mf'_y(y_0)(-M+f'_y(y_0)) g(y_0) -Mf''_{yy}(y_0)(g(y_0),g(y_0))\big)+h^4\Big(b_3(a_{31}+a_{32})a_{32}a_{21}\cr\noalign{\vskip2truemm}
&+b_4(a_{41}+a_{42}+a_{43})(a_{42}a_{21}+a_{43}(a_{31}+a_{32}))\Big)f''_{yy}(y_0)((-M+f'_y(y_0))g(y_0),g(y_0))\cr\noalign{\vskip2truemm}
&+h^4b_4a_{43}a_{32}a_{21} f'_y(y_0)(-M+f'_y(y_0))(-M+f'_y(y_0))g(y_0)+\frac{h^4}{2!}\big(b_3a_{32}a_{21}^2+b_4a_{42}a_{21}^2\cr\noalign{\vskip2truemm}
&+b_4a_{43}(a_{31}^2+2a_{31}a_{32}+a_{32}^2)\big)f'_y(y_0)f''_{yy}(y_0)\big(g(y_0),g(y_0)\big)+\frac{h^4}{3!}\Big(b_2a_{21}^3+b_3\big(a_{31}^3+3(a_{31}^2a_{32}
\cr\noalign{\vskip2truemm}
&+a_{31}a_{32}^2) +a_{32}^3\big)+b_4\big(a_{41}^3+a_{42}^3+a_{43}^3+3(a_{41}(a_{42}^2+a_{43}^2)+a_{42}(a_{41}^2+a_{43}^2)+a_{43}(a_{41}^2+a_{42}^2))\cr\noalign{\vskip2truemm}
&+6a_{41}a_{42}a_{43}\big)\Big)f'''_{yyy}(y_0)\big(g(y_0),g(y_0),g(y_0)\big)+\mathcal{O}{(h^5)}.
\end{aligned}
\end{equation*}}
 Under the assumption $c_i=\sum\limits_{j=1}^{i-1}a_{ij}$ for  $i=1,\ldots,4$, by comparing the Taylor series  of the exact solution $y(t_0+h)$,  we can verify that the explicit MVERK method \label{equ4} with coefficients  satisfying the conditions \eqref{order condition}, has order four. This completes the proof.
 \end{proof}

 \begin{rem}
Third-order  MVERK methods appearing in \cite{Wang2022}  use the Jacobian matrix of $f(y)$ with respect to $y$,  and the fourth-order MVERK method \eqref{MVERK4} uses the Jacobian matrix and the  Hessian matrix  of $f(y)$ with respect to $y$  at each step, however, as we known that this idea for first-order differential problems is no by means new (see, e.g., \cite{Abhulimen2014,Cash1981Siam,Enright1974siam}).
 It can be predicted that 
  high-order MVERK methods need to  use high-order derivative of of $f(y)$ with respect to $y$, hence our study does not involve high-order $(p\geq5)$  MVERK methods.
 \end{rem}
    
   The choice is  that  $a_{21}=\frac{1}{2}$, $a_{32}=\frac{1}{2}$,  $a_{43}=1$ and $b_1=\frac{1}{6}$, $b_2=\frac{2}{6}$, $b_3=\frac{2}{6}$, $b_4=\frac{1}{6}$, which  satisfies the order conditions (\ref{order condition}). Thus we obtain the  fourth-order explicit MVERK method with four stages as follows:
\begin{equation}\label{MVERK41}
\left\{
\begin{array}{l}
\displaystyle Y_1=y_0, \cr\noalign{\vskip3truemm}
\displaystyle Y_2=y_0+\frac{h}{2}(-MY_1+f(Y_1)),\cr\noalign{\vskip3truemm}
\displaystyle Y_3=y_0+\frac{h}{2}(-MY_2+f(Y_2)), \cr\noalign{\vskip3truemm}
\displaystyle Y_4=y_0+h(-MY_3+f(Y_3)), \cr\noalign{\vskip3truemm}
\displaystyle y_1=e^{-hM}y_0+\frac{h}{6}(f(Y_1)+2f(Y_2)+2f(Y_3)+f(Y_4))+w_4(z),
\end{array}
\right.
\end{equation}
the method (\ref{MVERK41}) can also  be expressed in the Butcher tableau
\begin{equation}\label{Tableau-MVERK41}
\begin{aligned}
\begin{tabular}{c|c|cccc}
$0$ &$I$&$0$&&\\
$\raisebox{-1.3ex}[1.0pt]{$\frac{1}{2}$}$& $\raisebox{-1.3ex}[1.0pt]{$I$}$ &$\raisebox{-1.3ex}[1.0pt]{$\frac{1}{2}$}$ &$\raisebox{-1.3ex}[1.0pt]{$0$}$&\\
 $\raisebox{-1.3ex}[1.0pt]{$\frac{1}{2}$}$ &$\raisebox{-1.3ex}[1.0pt]{$I$}$&  $\raisebox{-1.3ex}[1.0pt]{$0$}$& $\raisebox{-1.3ex}[1.0pt]{$\frac{1}{2}$}$& $\raisebox{-1.3ex}[1.0pt]{$0$}$&\\
 $\raisebox{-1.3ex}[1.0pt]{$1$}$ &$\raisebox{-1.3ex}[1.0pt]{$I$}$&  $\raisebox{-1.3ex}[1.0pt]{$0$}$& $\raisebox{-1.3ex}[1.0pt]{$0$}$& $\raisebox{-1.3ex}[1.0pt]{$1$}$& $\raisebox{-1.3ex}[1.0pt]{$0$}$\\[5pt]
 \hline
 $\raisebox{-1.3ex}[1.0pt]{$e^{-hM}$}$ &$\raisebox{-1.3ex}[1.0pt]{$w_4(z)$}$&$\raisebox{-1.3ex}[1.0pt]{$\frac{1}{6}$}$&\raisebox{-1.3ex}[1.0pt]{$\frac{2}{6}$} &  $\raisebox{-1.3ex}[1.0pt]{$\frac{2}{6}$}$&$\raisebox{-1.3ex}[1.0pt]{$\frac{1}{6}$}$\\
\end{tabular}
\end{aligned}
\end{equation}
The fourth-order explicit MVERK method (\ref{MVERK41}) reduces to the `classical Runge--Kutta method' when $M \rightarrow \mathbf{0}$,  which is especially notable.
 We here mention some fourth-order explicit ERK methods in the literature. Hochbruck et al. \cite{Hochbruck2005a}  proposed the five stages explicit ERK method of  order  four  which can be indicated by  the Butcher tableau
\begin{equation}\label{Tableau-ERK51}
\begin{aligned}
\begin{tabular}{c|ccccc}
$0$&&&\\
$\frac{1}{2}$&$\frac{1}{2}\varphi_{1,2}$&&\\
$\raisebox{-1.3ex}[1.0pt]{$\frac{1}{2}$}$&$\raisebox{-1.3ex}[1.0pt]{$\frac{1}{2}\varphi_{1,3}-\varphi_{2,3}$}$ &$\raisebox{-1.3ex}[1.0pt]{$\varphi_{2,3}$}$&\\
 $\raisebox{-1.3ex}[1.0pt]{$1$}$&  $\raisebox{-1.3ex}[1.0pt]{$\varphi_{1,4}-2\varphi_{2,4}$}$& $\raisebox{-1.3ex}[1.0pt]{$\varphi_{2,4}$}$& $\raisebox{-1.3ex}[1.0pt]{$\varphi_{2,4}$}$&\\
 $\raisebox{-1.3ex}[1.0pt]{$\frac{1}{2}$}$&  $\raisebox{-1.3ex}[1.0pt]{$\frac{1}{2}\varphi_{1,5}-2a_{5,2}-a_{5,4}$}$& $\raisebox{-1.3ex}[1.0pt]{$a_{5,2}$}$& $\raisebox{-1.3ex}[1.0pt]{$a_{5,2}$}$& $\raisebox{-1.3ex}[1.0pt]{$\frac{1}{4}\varphi_{2,5}-a_{5,2}$}$\\[8pt]
 \hline
&$\raisebox{-1.3ex}[1.0pt]{$\varphi_1-3\varphi_2+4\varphi_3$}$&\raisebox{-1.3ex}[1.0pt]{$0$} &  $\raisebox{-1.3ex}[1.0pt]{$0$}$&$\raisebox{-1.3ex}[1.0pt]{$-\varphi_2+4\varphi_3$}$&$\raisebox{-1.3ex}[1.0pt]{$4\varphi_2-8\varphi_3$}$\\
\end{tabular}
\end{aligned}
\end{equation}
with
$$a_{5,2}=\frac{1}{2}\varphi_{2,5}-\varphi_{3,4}+\frac{1}{4}\varphi_{2,4}-\frac{1}{2}\varphi_{3,5},$$
and Krogstad \cite{Krogstad2005} presented the fourth-order  ERK method with four stages which can be denoted by  the Butcher tableau
\begin{equation}\label{Tableau-ERK41}
\begin{aligned}
\begin{tabular}{c|cccc}
$0$&&\\
$\frac{1}{2}$&$\frac{1}{2}\varphi_{1,2}$&\\
$\raisebox{-1.3ex}[1.0pt]{$\frac{1}{2}$}$&$\raisebox{-1.3ex}[1.0pt]{$\frac{1}{2}\varphi_{1,3}-\varphi_{2,3}$}$ &$\raisebox{-1.3ex}[1.0pt]{$\varphi_{2,3}$}$&\\
 $\raisebox{-1.3ex}[1.0pt]{$1$}$&  $\raisebox{-1.3ex}[1.0pt]{$\varphi_{1,4}-2\varphi_{2,4}$}$& $\raisebox{-1.3ex}[1.0pt]{$0$}$& $\raisebox{-1.3ex}[1.0pt]{$2\varphi_{2,4}$}$&\\[6pt]
 \hline
&$\raisebox{-1.3ex}[1.0pt]{$\varphi_1-3\varphi_2+4\varphi_3$}$&\raisebox{-1.3ex}[1.0pt]{$2\varphi_2-4\varphi_3$} &  $\raisebox{-1.3ex}[1.0pt]{$2\varphi_2-4\varphi_3$}$&$\raisebox{-1.3ex}[1.0pt]{$-\varphi_2+4\varphi_3$}$\\
\end{tabular}
\end{aligned}
\end{equation}
where
\begin{equation}
\varphi_{i,j}=\varphi_i(-c_jhM)=\int_{0}^1 e^{-(1-\tau)c_jhM}\frac{\tau^{i-1}}{(i-1)!}d\tau.
\end{equation}
 As claimed by Hochbruck et al. \cite{Hochbruck2005a}, the methods  \eqref{Tableau-ERK51}   and
 \eqref{Tableau-ERK41} do not satisfy the stiff-order conditions strictly, but to a weak form.
  The coefficients of (\ref{Tableau-ERK51}) and (\ref{Tableau-ERK41}) are  matrix exponentials, normally, it is needed to recalculate the coefficients of them at each time step once we consider the variable stepsize technique in practice. The coefficients $a_{ij}$, $b_i$, for $i,j=1, \ldots, s$ of  MVERK methods are real constants, which can greatly reduce the computational cost of matrix exponentials. 

 The another choice of $a_{21}=\frac{1}{3}$, $a_{31}=-\frac{1}{3}$, $a_{32}=1$, $a_{41}=1$, $a_{42}=-1$, $a_{43}=1$ and $b_1=\frac{1}{8}$, $b_2=\frac{3}{8}$, $b_3=\frac{3}{8}$, $b_4=\frac{1}{8}$, which also satisfies the order conditions (\ref{order condition}). Therefore we get another fourth-order expkicit MVERK method with four stages as follows:
\begin{equation}\label{MVERK42}
\left\{
\begin{array}{l}
\displaystyle Y_1=y_0, \cr\noalign{\vskip3.5truemm}
\displaystyle Y_2=y_0+\frac{h}{3}(-MY_1+f(Y_1)),\cr\noalign{\vskip3.5truemm}
\displaystyle Y_3=y_0-\frac{h}{3}(-MY_1+f(Y_1))+h(-MY_2+f(Y_2)), \cr\noalign{\vskip3.5truemm}
\displaystyle Y_4=y_0+h(-MY_1+f(Y_1))-h(-MY_2+f(Y_2))+h(-MY_3+f(Y_3)), \cr\noalign{\vskip3.5truemm} \displaystyle y_{1}=e^{-hM}y_0+\frac{h}{8}(f(Y_1)+3f(Y_2)+3f(Y_3)+f(Y_4))+w_4(z).
\end{array}
\right.
\end{equation}
The method (\ref{MVERK42}) can be expressed in the  Butcher tableau
\begin{equation}\label{Tableau-MVERK42}
\begin{aligned}
\begin{tabular}{c|c|cccc}
$0$ &$I$&$0$&&\\
$\raisebox{-1.3ex}[1.0pt]{$\frac{1}{3}$}$& $\raisebox{-1.3ex}[1.0pt]{$I$}$ &$\raisebox{-1.3ex}[1.0pt]{$\frac{1}{3}$}$ &$\raisebox{-1.3ex}[1.0pt]{$0$}$&\\
 $\raisebox{-1.3ex}[1.0pt]{$\frac{2}{3}$}$ &$\raisebox{-1.3ex}[1.0pt]{$I$}$&  $\raisebox{-1.3ex}[1.0pt]{-$\frac{1}{3}$}$& $\raisebox{-1.3ex}[1.0pt]{$1$}$& $\raisebox{-1.3ex}[1.0pt]{$0$}$&\\
 $\raisebox{-1.3ex}[1.0pt]{$1$}$ &$\raisebox{-1.3ex}[1.0pt]{$I$}$&  $\raisebox{-1.3ex}[1.0pt]{$1$}$& $\raisebox{-1.3ex}[1.0pt]{-$1$}$& $\raisebox{-1.3ex}[1.0pt]{$1$}$& $\raisebox{-1.3ex}[1.0pt]{$0$}$\\[5pt]
 \hline
 $\raisebox{-1.3ex}[1.0pt]{$e^{-hM}$}$ &$\raisebox{-1.3ex}[1.0pt]{$w_4(z)$}$&$\raisebox{-1.3ex}[1.0pt]{$\frac{1}{8}$}$&\raisebox{-1.3ex}[1.0pt]{$\frac{3}{8}$} &  $\raisebox{-1.3ex}[1.0pt]{$\frac{3}{8}$}$&$\raisebox{-1.3ex}[1.0pt]{$\frac{1}{8}$}$\\
\end{tabular}
\end{aligned}
\end{equation}
Likewise, when $M \rightarrow 0$, the fourth-order explicit MVERK method (\ref{MVERK42}) reduces to the Kutta's fourth-order method, or `3/8-Rule'.

\section{A simplified version of fourth-order explicit ERK methods}\label{sec3}
Differently from MVERK methods, the internal stages and the update of the simplified version of ERK methods simultaneously preserve  some properties of matrix exponentials.
\begin{mydef}\label{definition2}(\cite{Wang2022})
An $s$-stage simplified version of   exponential
{Runge--Kutta} (SVERK) method for the numerical integration
(\ref{equ1}) is defined as
\begin{equation}\label{SVERK}
\left\{
\begin{array}{l}
\displaystyle Y_i=e^{-\bar{c}_ihM}y_0+h\sum\limits_{j=1}^s\bar{a}_{ij}f(Y_j), \
i=1,\ldots,s,\cr\noalign{\vskip3truemm}
\displaystyle y_{1}=e^{-hM}y_0+h\sum\limits_{i=1}^s\bar{b}_{i}f(Y_i)+\bar {w}_s(z),
\end{array}
\right.
\end{equation}
{where} $\bar{a}_{ij}$, $\bar{b}_i$
 are real constants for ${i,j=1,\ldots s}$, $Y_i\approx y(t_0+\bar{c}_ih) $ for $i=1,\ldots,s$, $\bar{w}_s(z)$ is related to $hM$, and $\bar{w}_s(z)\rightarrow0$ when $M\rightarrow \mathbf{0}$.
\end{mydef}

 Similarly to $w_s(z)$, here $\bar{w}_s(z)$  is independent of the matrix exponential, and  SVERK methods with the same order share the same $\bar{w}_s(z)$.  Once we consider the order of SVERK methods which satisfies $p\geq2$, then the $\bar{w}_s(z)$ is  related to the term $f(\cdot)$ and initial value $y_0$.  In particular, when we consider the first-order SVERK method, then $\bar{w}_s(z)=0$. It is a fact that   the $\bar{w}_s(z)$ of SVERK methods is different from the $w_s(z)$ of MVERK  methods when $p \geq 3$, and we  have shown the difference of  $\bar{w}_3(z)$ and $w_3(z)$ in  \cite{Wang2022}. In what follows, we will present   the difference between $\bar{w}_4(z)$ and $w_4(z)$.

The method (\ref{SVERK}) can be displayed by the following Butcher Tableau

\begin{equation}\label{Tableau-SVERK4}
\begin{aligned} &\quad\quad\begin{tabular}{c|c|c}
 ${\bar{c}}$&$\mathbf{e}^{-{\bar{c}}hM}$&${\bar{A}}$ \\
 \hline
  $\raisebox{-1.3ex}[1.0pt]{$e^{-hM}$}$ & $\raisebox{-1.3ex}[1.0pt]{$\bar{w}_s(z)$}$&$\raisebox{-1.3ex}[1.0pt]{$\bar{b}^{\mathrm{T}}$}$  \\
\end{tabular}
~=
\begin{tabular}{c|c|ccc}
 $\bar{c}_1$&$e^{-\bar{c}_1hM}$&$\bar{a}_{11}$&$\cdots$&$\bar{a}_{1s}$\\
$\vdots$& $\vdots$ & $\vdots$&$\vdots$&$\vdots$\\
 $\bar{c}_s$ &$e^{-\bar{c}_shM}$&  $\bar{a}_{s1}$& $\cdots$& $\bar{a}_{ss}$\\
 \hline
 $\raisebox{-1.3ex}[1.0pt]{$e^{-hM}$}$&$\raisebox{-1.3ex}[1.0pt]{$\bar{w}_s(z)$}$&$\raisebox{-1.3ex}[1.0pt]{$\bar{b}_1$}$&\raisebox{-1.3ex}[1.0pt]{$\cdots$} &  $\raisebox{-1.3ex}[1.0pt]{$\bar{b}_s$}$\\
\end{tabular}
\end{aligned}
\end{equation}
with the suitable assumption $\bar{c}_i=\sum\limits_{j=1}^s \bar{a}_{ij}, \ i=1,\ldots,s$. From (\ref{Tableau-SVERK4}),  the coefficients of the SVERK method are real constants, which are different from  standard ERK methods (see, e.g., \cite{Hochbruck2005a,Hochbruck2005b,Hochbruck2010}). 
Obviously,  SVERK methods can integrate the first-order homogeneous linear system  exactly, so do  MVERK methods. 
 
 The following theorem will show the order conditions of fourth-order explicit 
SVERK methods with four stages.

\begin{mytheo}
If  the coefficients of the four-stage explicit SVERK method
 \begin{equation}\label{SVERK44}
\left\{
\begin{array}{l}
\displaystyle Y_1=y_0, \cr\noalign{\vskip3truemm}
\displaystyle Y_2=e^{-\bar{c}_2hM}y_0+h\bar{a}_{21}f(Y_1),\cr\noalign{\vskip3truemm}
\displaystyle Y_3=e^{-\bar{c}_3hM}y_0+h\big(\bar{a}_{31}f(Y_1)+\bar{a}_{32}f(Y_2)\big), \cr\noalign{\vskip3truemm}
\displaystyle Y_4=e^{-\bar{c}_4hM}y_0+h\big(\bar{a}_{41}f(Y_1)+\bar{a}_{42}f(Y_2)+\bar{a}_{43}f(Y_3)\big), \cr\noalign{\vskip3truemm}
\displaystyle y_1=e^{-hM}y_0+h\big(\bar{b}_1f(Y_1)+\bar{b}_2f(Y_2)+\bar{b}_3f(Y_3)+\bar{b}_4f(Y_4)\big)+\bar{w}_4(z),
\end{array}
\right.
\end{equation}
 where
 \small{
 \begin{equation*}
 \begin{aligned}
 \bar{w}_4(z)&=-\frac{h^2}{2!}Mf(y_0)+\frac{h^3}{3!}((M-f'_y(y_0))Mf(y_0)-Mf'_y(y_0)g(y_0))+\frac{h^4}{4!}\big((-M+f'_y(y_0))\cr\noalign{\vskip2truemm}
& \quad \cdot M^2f(y_0)+M^2f'_y(y_0)g(y_0)-Mf''_{yy}(y_0)\big(g(y_0),g(y_0)\big)-Mf'_y(y_0)(-M+f'_y(y_0)) g(y_0)\cr\noalign{\vskip2truemm}
&\quad -f'_y(y_0)Mf'_y(y_0)g(y_0)-f'_y(y_0)f'_y(y_0)Mf(y_0)+3f''_{yy}(y_0)(-Mf(y_0),g(y_0))\big),
 \end{aligned}
 \end{equation*}}
and $g(y_0)=-My_0+f(y_0)$, satisfy the order conditions \eqref{order condition}, then the explicit SVERK method has order four.
\end{mytheo}
\begin{proof}
Under the  assumptions  $y(t_0)=y_0$ and $\bar{c}_i=\sum\limits_{j=1}^{i-1}\bar{a}_{ij}$ for $i=1,2,3,4$, the Taylor series of the numerical solution $y_1$ is given by
\begin{equation*}
\begin{aligned}
y_1=&\big(I-hM+\frac{(hM)^2}{2!}-\frac{(hM)^3}{3!}+\frac{(hM)^4}{4!}+\mathcal{O}{(h^5)}\big)y_0+h\Big(\bar{b}_1f(y_0)+
\bar{b}_2f\big((I-\bar{c}_2hM\cr\noalign{\vskip3truemm}
&+\frac{(\bar{c}_2hM)^2}{2!}-\frac{(\bar{c}_2hM)^3}{3!}+\mathcal{O}{(h^4)})y_0+h\bar{a}_{21}f(y_0)\big)+\bar{b}_3f\big((I-\bar{c}_3hM+\frac{(\bar{c}_3hM)^2}{2!}\cr\noalign{\vskip3truemm}
&-\frac{(\bar{c}_3hM)^3}{3!}+\mathcal{O}{(h^4)})y_0+h\bar{a}_{31}f(y_0)+h\bar{a}_{32}f(Y2)\big)+\bar{b}_4f\big((I-\bar{c}_4hM+\frac{(\bar{c}_4hM)^2}{2!}\cr\noalign{\vskip3truemm}
&-\frac{(\bar{c}_4hM)^3}{3!}+\mathcal{O}{(h^4)})y_0+h\bar{a}_{41}f(y_0)+h\bar{a}_{42}f(Y2)+h\bar{a}_{43}f(Y3)\big)\Big)+\bar{w}_4(z).
\end{aligned}
\end{equation*}
Inserting  the internal stages $Y_2, \ Y_3,\ Y_4$ into the numerical solution $y_1$, and expanding $f(Y_2),\ f(Y_3),\ f(Y_4)$  in a Taylor series with  respect to $y_0$ leads to
\begin{equation*}
\begin{aligned}
y_1=&\big(I-hM+\frac{(hM)^2}{2!}-\frac{(hM)^3}{3!}+\frac{(hM)^4}{4!}+\mathcal{O}{(h^5)}\big)y_0+h\bar{b}_1f(y_0)+h\bar{b}_2\Big(f(y_0)+hf_y(y_0)\cr\noalign{\vskip2truemm}
&\cdot (-\bar{c}_2My_0+\bar{a}_{21}f(y_0)+\frac{h(\bar{c}_2M)^2y_0}{2!}-\frac{h^2(\bar{c}_2M)^3y_0}{3!})+\frac{h^2}{2}f''_{yy}(y_0)\big(-\bar{c}_2My_0+\bar{a}_{21}f(y_0)\cr\noalign{\vskip2truemm}
&+\frac{h(\bar{c}_2M)^2y_0}{2!},-\bar{c}_2My_0+\bar{a}_{21}f(y_0)+\frac{h(\bar{c}_2M)^2y_0}{2!}\big)+\frac{h^3}{3}f'''_{yyy}(y_0)\big(-\bar{c}_2My_0+\bar{a}_{21}f(y_0),\cr\noalign{\vskip2truemm}
&-\bar{c}_2My_0+\bar{a}_{21}f(y_0),-\bar{c}_2My_0+\bar{a}_{21}f(y_0)\big)\Big)+h\bar{b}_3\Big(f(y_0)+hf_y(y_0) \big(-\bar{c}_3My_0+\bar{a}_{31}f(y_0)\cr\noalign{\vskip2truemm}
&+\bar{a}_{32}f(Y2)+\frac{h(\bar{c}_3M)^2y_0}{2!} -\frac{h^2(\bar{c}_3M)^3y_0}{3!}\big)+\frac{h^2}{2!}f''_{yy}(y_0)\big(-\bar{c}_3My_0+\bar{a}_{31}f(y_0) +\bar{a}_{32}f(Y2)\cr\noalign{\vskip2truemm}
& +\frac{h(\bar{c}_3M)^2y_0}{2!},-\bar{c}_3My_0 +\bar{a}_{31}+\bar{a}_{32}f(Y2)+\frac{h(\bar{c}_3M)^2y_0}{2!}\big)+\frac{h^3}{3!}f'''_{yyy}(y_0)\big(-\bar{c}_3My_0+\bar{a}_{31}\cr\noalign{\vskip2truemm}
&\cdot f(y_0)+\bar{a}_{32}f(Y2),-\bar{c}_3My_0+\bar{a}_{31}f(y_0)+\bar{a}_{32}f(Y2),-\bar{c}_3My_0+\bar{a}_{31}f(y_0)+\bar{a}_{32}f(Y2)\big)\Big)\cr\noalign{\vskip2truemm}
&+h\bar{b}_4\Big(f(y_0)+hf_y(y_0) \big(-\bar{c}_4My_0+\bar{a}_{41}f(y_0)+\bar{a}_{42}f(Y2)+\bar{a}_{43}f(Y_3)+\frac{h(\bar{c}_4M)^2y_0}{2!} \cr\noalign{\vskip2truemm}
&-\frac{h^2(\bar{c}_4M)^3y_0}{3!}\big)+\frac{h^2}{2!}f''_{yy}(y_0)\big(-\bar{c}_4My_0+\bar{a}_{41}f(y_0)+\bar{a}_{42}f(Y2)+\bar{a}_{43}f(Y_3) +\frac{h(\bar{c}_4M)^2y_0}{2!},\cr\noalign{\vskip2truemm}
&-\bar{c}_4My_0 +\bar{a}_{41}f(y_0)+\bar{a}_{42}f(Y2)+\bar{a}_{43}f(Y_3)+\frac{h(\bar{c}_4M)^2y_0}{2!}\big)\cr\noalign{\vskip2truemm}
&+\frac{h^3}{3!}f'''_{yyy}(y_0)\big(-\bar{c}_4My_0+\bar{a}_{41}f(y_0)+\bar{a}_{42}f(Y2)+\bar{a}_{43}f(Y_3),-\bar{c}_4My_0+\bar{a}_{41}f(y_0)+\bar{a}_{43}f(Y_3),\cr\noalign{\vskip2truemm}
&-\bar{c}_4My_0+\bar{a}_{41}f(y_0)+\bar{a}_{42}f(Y2)+\bar{a}_{43}f(Y_3)\big)\Big)+\bar{w}_4(z).
\end{aligned}
\end{equation*}
To sum up, we have
\begin{equation*}
\begin{aligned}
y_1=&y_0-hMy_0+h(\bar{b}_1+\bar{b}_2+\bar{b}_3+\bar{b}_4)f(y_0)-\frac{h^2Mg(y_0)}{2!}+h^2(\bar{b}_2\bar{c}_{2}+\bar{b}_{3}\bar{c}_{3}+\bar{b}_{4}\bar{c}_{4}) f'_y(y_0)g(y_0) \cr\noalign{\vskip3truemm}
&-\frac{h^3}{3!}f'_y(y_0)Mf(y_0)+\frac{h^3}{3!}M(M-f'_y(y_0))g(y_0)+\frac{h^3}{2!}(\bar{b}_2\bar{a}_{21}^2+\bar{b}_3\bar{c}_{3}^2+\bar{b}_4\bar{c}_{4}^2) f''_{yy}(y_0)(g(y_0),\cr\noalign{\vskip3truemm}
&g(y_0)) +h^3(\bar{b}_3\bar{a}_{32}\bar{a}_{21}+\bar{b}_4\bar{a}_{42}\bar{a}_{21}+\bar{b}_4\bar{a}_{43}\bar{c}_{3})f'_y(y_0)(M^2y_0+f'_y(y_0)g(y_0))\cr\noalign{\vskip3truemm}
&+\frac{h^4}{4!}\Big(-M^3g(y_0)+M^2f'_y(y_0)g(y_0)-Mf'_y(y_0)(-M+f'_y(y_0))g(y_0)-Mf''_{yy}(y_0)(g(y_0),\cr\noalign{\vskip3truemm}
&g(y_0))-f'_y(y_0)Mf'_y(y_0)g(y_0)-f'_y(y_0)f'_y(y_0)Mf(y_0)+3f''_{yy}(y_0)\big(-Mf(y_0),g(y_0)\big)\Big)\cr\noalign{\vskip3truemm}
&+\frac{h^4}{3!}\big(\bar{b}_2\bar{c}_{2}^3
\bar{b}_3\bar{c}_{3}^3+\bar{b}_{4}\bar{c}_{4}^3\big)f'''_{yyy}(y_0)\big(g(y_0),g(y_0),g(y_0)\big)
+h^4\Big(\bar{b}_3\bar{c}_{3}\bar{a}_{32}\bar{a}_{21}+\bar{b}_4\bar{c}_{4}(\bar{a}_{42}\bar{a}_{21}\cr\noalign{\vskip3truemm}
&+\bar{a}_{43}\bar{c}_{3})\Big)f''_{yy}(y_0)\big(M^2 y_0+f'_y(y_0)g(y_0),g(y_0)\big)+\frac{h^4}{2!}\big(\bar{b}_3\bar{a}_{32}\bar{a}_{21}^2+\bar{b}_4\bar{a}_{42}\bar{a}_{21}^2+\bar{b}_4\bar{a}_{43}\bar{c}_{3}^2\big) \cr\noalign{\vskip3truemm}
&\cdot  f'_y(y_0)f''_{yy}(y_0)\big(g(y_0),g(y_0)\big)+h^4\bar{b}_4\bar{a}_{43}\bar{a}_{32}\bar{a}_{21}f'_y(y_0)\big((M^2+f'_y(y_0)^2)g(y_0)\cr\noalign{\vskip3truemm}
& +f'_y(y_0)^2M^2y_0\big)+\mathcal{O}(h^5).
\end{aligned}
\end{equation*}
Hence, the coefficients of the four-stage explicit SVERK method \eqref{SVERK44} satisfying  the order conditions \eqref{order condition}, has order four. This completes the proof.
 \end{proof}
 
\begin{rem}
Compared with fourth-order explicit MVERK methods, the internal stages of fourth-order explicit SVERK methods use matrix exponentials of $hM$ to some extent, and the update  also uses the Jacobian matrix and the Hessian matrix of $f(y)$  with respect to $y$ at each step. Once we consider the high-order SVERK method, and the high-order derivative of $f(y)$ with respect to $y$ is needed. Hence, this paper does not involve  high-order $(p\geq5)$ SVERK methods, and it can be observed that the $\bar{w}_{4}(z)$ of \eqref{SVERK44}   is more complicated than the $w_{4}(z)$ of \eqref{MVERK4}, and our numerical experiments will support this point by comparing their consumed CPU times (in seconds).
\end{rem}

 Taking $a_{21}=\frac{1}{2}, a_{32}=\frac{1}{2}, a_{43}=1$ and $b_1=\frac{1}{6}, b_2=\frac{2}{6}, b_3=\frac{2}{6}, b_4=\frac{1}{6}$, we  then obtain the fourth-order explicit SVERK method with four stages:
 \begin{equation}\label{SVERK41}
\left\{
\begin{array}{l}
Y_1=y_0, \cr\noalign{\vskip3.5truemm}
 Y_2=e^{-\frac{1}{2}hM}y_0+\frac{h}{2}f(Y_1),\cr\noalign{\vskip3.5truemm}
Y_3=e^{-\frac{1}{2}hM}y_0+\frac{h}{2}f(Y_2), \cr\noalign{\vskip3.5truemm}
 \displaystyle Y_4=e^{-hM}y_0+hf(Y_3), \cr\noalign{\vskip3.5truemm}
  \displaystyle y_1=e^{-hM}y_0+\frac{h}{6}(f(Y_1)+2f(Y_2)+2f(Y_3)+f(Y_4))+\bar{w}_4(z),
\end{array}
\right.
\end{equation}
with $g(y_0)=-My_0+f(y_0).$
The method (\ref{SVERK41}) can be indicated by the Butcher tableau
\begin{equation}\label{Tableau-SVERK41}
\begin{aligned}
\begin{tabular}{c|c|cccc}
$0$ &$I$&$0$&&\\
$\raisebox{-1.3ex}[1.0pt]{$\frac{1}{2}$}$& $\raisebox{-1.3ex}[1.0pt]{$e^{-\frac{1}{2}hM}$}$ &$\raisebox{-1.3ex}[1.0pt]{$\frac{1}{2}$}$ &$\raisebox{-1.3ex}[1.0pt]{$0$}$&\\
 $\raisebox{-1.3ex}[1.0pt]{$\frac{1}{2}$}$ &$\raisebox{-1.3ex}[1.0pt]{$e^{-\frac{1}{2}hM}$}$&  $\raisebox{-1.3ex}[1.0pt]{$0$}$& $\raisebox{-1.3ex}[1.0pt]{$\frac{1}{2}$}$& $\raisebox{-1.3ex}[1.0pt]{$0$}$&\\
 $\raisebox{-1.3ex}[1.0pt]{$1$}$ &$\raisebox{-1.3ex}[1.0pt]{$e^{-hM}$}$&  $\raisebox{-1.3ex}[1.0pt]{$0$}$& $\raisebox{-1.3ex}[1.0pt]{$0$}$& $\raisebox{-1.3ex}[1.0pt]{$1$}$& $\raisebox{-1.3ex}[1.0pt]{$0$}$\\[5pt]
 \hline
 $\raisebox{-1.3ex}[1.0pt]{$e^{-hM}$}$ &$\raisebox{-1.3ex}[1.0pt]{$\bar{w}_4(z)$}$&$\raisebox{-1.3ex}[1.0pt]{$\frac{1}{6}$}$&\raisebox{-1.3ex}[1.0pt]{$\frac{2}{6}$} &  $\raisebox{-1.3ex}[1.0pt]{$\frac{2}{6}$}$&$\raisebox{-1.3ex}[1.0pt]{$\frac{1}{6}$}$\\
\end{tabular}
\end{aligned}
\end{equation}

 Another option is that  $a_{21}=\frac{1}{3}$, $a_{31}=-\frac{1}{3}$, $a_{32}=1$, $a_{41}=1$, $a_{42}=-1$, $a_{43}=1$, $b_1=\frac{1}{8}$, $b_2=\frac{3}{8}$, $b_3=\frac{3}{8}$  and  $ b_4=\frac{1}{8}$. Thus, we have the following fourth-order explicit SVERK method with four stages:
\begin{equation}\label{SVERK42}
\left\{
\begin{array}{l}
 Y_1=y_0, \cr\noalign{\vskip3truemm}
  Y_2=e^{-\frac{1}{3}hM}y_0+\frac{h}{3}f(Y_1),\cr\noalign{\vskip3truemm}
Y_3=e^{-\frac{2}{3}hM}y_0-\frac{h}{3}f(Y_1)+hf(Y2), \cr\noalign{\vskip3truemm}
 Y_4=e^{-hM}y_0+hf(Y_1)-hf(Y_2)+hf(Y_3), \cr\noalign{\vskip3truemm} \displaystyle y_1=e^{-hM}y_0+\frac{h}{8}(f(Y_1)+3f(Y_2)+3f(Y_3)+f(Y_4))+\bar{w}_4(z),
\end{array}
\right.
\end{equation}
which can be represented by the Butcher tableau
\begin{equation}\label{Tableau-SVERK42}
\begin{aligned}
\begin{tabular}{c|c|cccc}
$0$ &$I$&$0$&&\\
$\raisebox{-1.3ex}[1.0pt]{$\frac{1}{3}$}$& $\raisebox{-1.3ex}[1.0pt]{$e^{-\frac{1}{3}hM}$}$ &$\raisebox{-1.3ex}[1.0pt]{$\frac{1}{3}$}$ &$\raisebox{-1.3ex}[1.0pt]{$0$}$&\\
 $\raisebox{-1.3ex}[1.0pt]{$\frac{2}{3}$}$ &$\raisebox{-1.3ex}[1.0pt]{$e^{-\frac{2}{3}hM}$}$&  $\raisebox{-1.3ex}[1.0pt]{-$\frac{1}{3}$}$& $\raisebox{-1.3ex}[1.0pt]{$1$}$& $\raisebox{-1.3ex}[1.0pt]{$0$}$&\\
 $\raisebox{-1.3ex}[1.0pt]{$1$}$ &$\raisebox{-1.3ex}[1.0pt]{$e^{-hM}$}$&  $\raisebox{-1.3ex}[1.0pt]{$1$}$& $\raisebox{-1.3ex}[1.0pt]{-$1$}$& $\raisebox{-1.3ex}[1.0pt]{$1$}$& $\raisebox{-1.3ex}[1.0pt]{$0$}$\\[5pt]
 \hline
 $\raisebox{-1.3ex}[1.0pt]{$e^{-hM}$}$ &$\raisebox{-1.3ex}[1.0pt]{$\bar{w}_4(z)$}$&$\raisebox{-1.3ex}[1.0pt]{$\frac{1}{8}$}$&\raisebox{-1.3ex}[1.0pt]{$\frac{3}{8}$} &  $\raisebox{-1.3ex}[1.0pt]{$\frac{3}{8}$}$&$\raisebox{-1.3ex}[1.0pt]{$\frac{1}{8}$}$\\
\end{tabular}
\end{aligned}
\end{equation}
Once $M \rightarrow 0$,  fourth-order explicit SVERK methods (\ref{SVERK41}) and (\ref{SVERK42}) with four stages reduce to  classical fourth-order explicit  RK methods.
\section{Convergence analysis}\label{sec4}

It is well known that the linear stability analysis of a RK method is to apply the method to a Dalquist equation \cite{Hairer2006}.  There is no doubt that we apply these exponential integrators to $y^{\prime}=i\lambda y,\ \lambda \in \mathbb R$, which are $A$-stable. Based on the work of  Buvoli and Minion \cite{Buvoli2022}, we apply these exponential methods to the partitioned Dalquist equation
\begin{equation}\label{test equation}
y^{\prime}=i\lambda_1y+i\lambda_2y, \qquad y(t_0)=y_0,\qquad \lambda_1, \lambda_2 \in \mathbb R.
\end{equation}
Solving the partitioned Dalquist equation \eqref{test equation} by a partitioned exponential integrator, and treating  the $i\lambda_1$ exponentially and the $i\lambda_2$ explicitly lead to
the  explicit scalar form
\begin{equation}
y_{n+1}=R(ik_1,ik_2)y_n, \quad \quad \ k_1=h\lambda_1, \ k_2=h\lambda_2.
\end{equation}
\begin{figure}[!htb]
\centering
\begin{tabular}[c]{cccc}%
  \subfigure[]{\includegraphics[width=6.6cm,height=5.5cm]{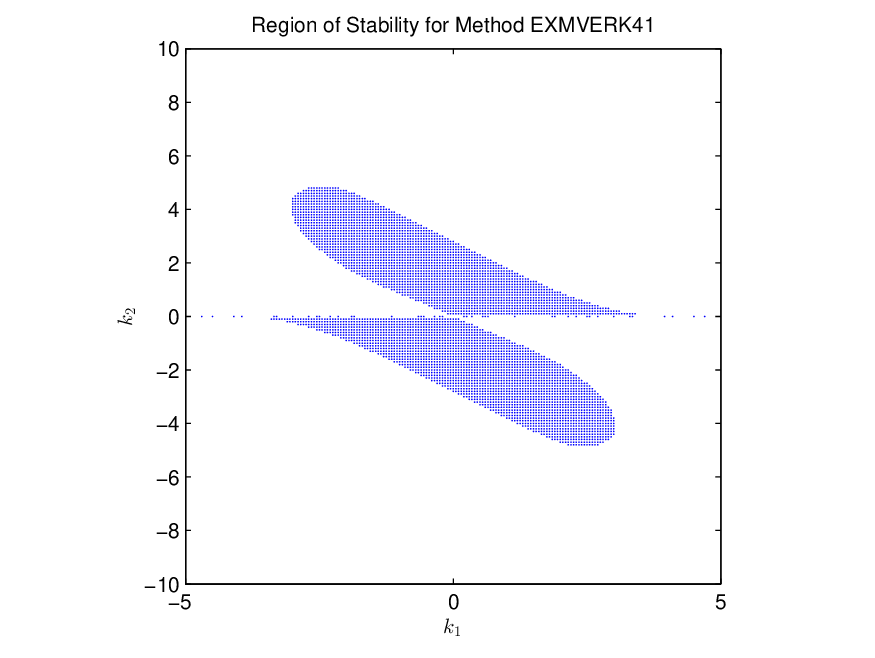}}
  \subfigure[]{\includegraphics[width=6.6cm,height=5.5cm]{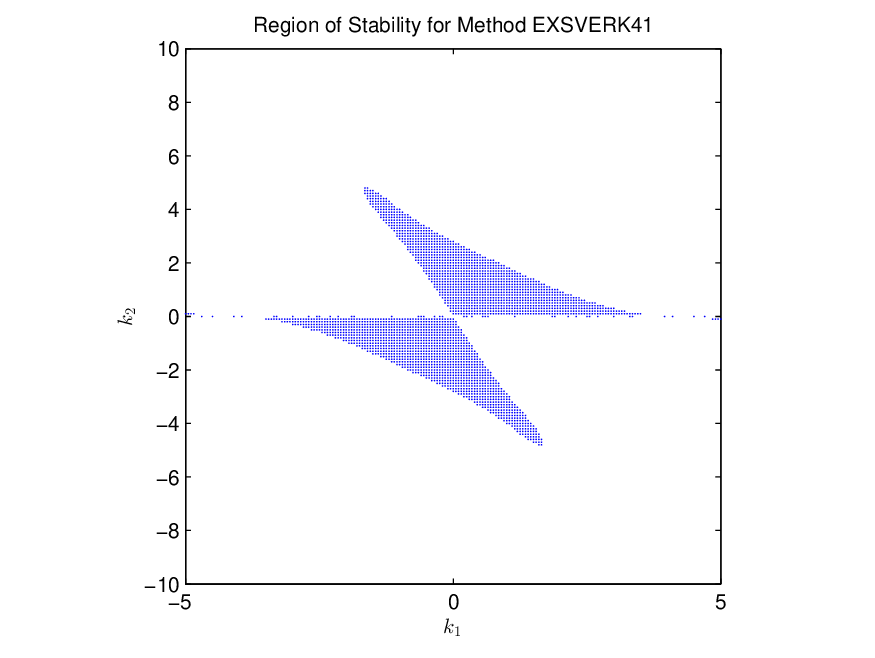}}
\end{tabular}
\caption{  {(a)}: Stability region for the fourth-order explicit MVERK (EXMVERK41) method \eqref{MVERK41} with four stages.  {(b)}: Stability region for the fourth-order explicit SVERK (EXSVERK41) method \eqref{SVERK41} with four stages.}\label{stability}
\end{figure}
 The stability regions of  fourth-order explicit ERK method \eqref{MVERK41}  and \eqref{SVERK41} are respectively depicted in Fig. \ref{stability} (a) and (b). The stability region of the standard  fourth-order   ERK method  \eqref{Tableau-ERK41} was shown in \cite{Buvoli2022},  and it can be observed that the stability regions of our methods are smaller than the stability region of  \eqref{Tableau-ERK41}.  
 
 In \cite{Wang2022}, the convergence of first-order explicit ERK methods was analyzed, however,  the convergence of high-order explicit ERK methods is not discussed.  The following theorem will present the  convergence of fourth-order explicit SVERK methods.

\begin{assum}\label{assum1}
Suppose that  $f$ satisfies  the local Lipschitz condition
in a strip along the exact solution $y$, i.e., there exists a real number $L>0$ for all $t\in[t_0,t_{end}]$, such that 
\begin{equation*}
\Vert f(\bar{y})-f(\hat{y})\Vert \leq L\Vert \bar{y}-\hat{y}\Vert,  \qquad \forall \hat{y}, \bar{y} \in T_{y},
\end{equation*}
where $T_{y}=\{\bar{y}\in V: min_{t} \|y(t)-\bar{y}\|\leq 1\}$ with the Euclidean norm $\|\cdot\|$.
\end{assum}
\begin{mytheo}
Under Assumption \ref{assum1}, if the  four-stage explicit
 SVERK method with $\bar{w}_4(z)$ satisfies the order conditions (\ref{order condition}),
  then the following error bound holds for $1\leq n \leq \frac{t_{end}-t_0}{h}$:
 $$\Vert y_n-y(t_n) \Vert \leq Ch^4,$$
  the constant $C$ is independent of $n$ and $h$, but  depends on $t_{end}$, $\Vert M \Vert$, $f^{(p)}(t_n)$ for $p=0,\ldots,4.$
\end{mytheo}
\begin{proof}
We rewrite the four-stage explicit SVERK method (\ref{SVERK44}) as
\begin{equation}\label{SVERK444}
\left\{
\begin{array}{l}
\displaystyle k_{n,1}=f(y_n), \cr\noalign{\vskip4truemm}
\displaystyle k_{n,2}=f(e^{-\bar{c}_2hM}y_n+h\bar{a}_{21}k_{n,1}),
\cr\noalign{\vskip4truemm}
\displaystyle k_{n,3}=f\big(e^{-\bar{c}_3hM}y_n+h(\bar{a}_{31}k_{n,1}+\bar{a}_{32}k_{n,2})\big)
, \cr\noalign{\vskip4truemm}
\displaystyle k_{n,4}=f\big(e^{-\bar{c}_4hM}y_n+h(\bar{a}_{41}k_{n,1}+\bar{a}_{42}k_{n,2}
+\bar{a}_{43}k_{n,3})\big)
, \cr\noalign{\vskip4truemm}
\displaystyle y_{n+1}=e^{-hM}y_n+h\big(\bar{b}_1k_{n,1}+\bar{b}_2k_{n,2}+\bar{b}_3k_{n,3}
+\bar{b}_4k_{n,4}\big)+\bar{w}_4(z),
\end{array}
\right.
\end{equation}
and expand  $y(t_n+h)$ into  a Taylor series
\begin{equation}
y(t_n+h)=y(t_n)+hy'(t_n)+\frac{h^2}{2!}y''(t_n)+\frac{h^3}{3!}y'''(t_n)
+\frac{h^4}{4!}y^{(4)}(t_n)+\mathcal{O}(h^5).
\end{equation}
Inserting the exact solution into the numerical scheme gives
  \begin{equation}
y(t_n+h)=e^{-hM}y(t_n)+h\big(\bar{b}_1\bar{k}_{n,1}+\bar{b}_2\bar{k}_{n,2}
+\bar{b}_3\bar{k}_{n,3}+\bar{b}_4\bar{k}_{n,4}\big)+\hat{w}_4(z)+\delta_{n+1},
\end{equation}
where
\begin{equation}\label{residual}
\begin{aligned}
\displaystyle \delta_{n+1}=&y(t_n)+hy'(t_n)+\frac{h^2}{2!}y''(t_n)+\frac{h^3}{3!}y'''(t_n)
+\frac{h^4}{4!}y^{(4)}(t_n)+\mathcal{O}(h^5)-e^{-hM}y(t_n)\cr\noalign{\vskip3truemm}
&-h\big(\bar{b}_1\bar{k}_{n,1}+\bar{b}_2\bar{k}_{n,2}
+\bar{b}_3\bar{k}_{n,3}+\bar{b}_4\bar{k}_{n,4}\big)-\hat{w}_4(z),
\end{aligned}
\end{equation}
and $\bar{k}_{n,1},\ \bar{k}_{n,2},\ \bar{k}_{n,3},\ \bar{k}_{n,4}$ and $\hat{w}_4(z)$ corresponding to $k_{n,,1},\ k_{n,2},\ k_{n,3},\ k_{n,4}$ and  $\bar{w}_4(z)$ satisfy $y(t_n)=y_n$, respectively.

If we denote $E_{n-1,i}=k_{n-1,i}-\bar{k}_{n-1,i}$ and $e_{n}=y_{n}-y(t_{n})$ , then
\begin{equation}\label{err1}
\begin{aligned}
\displaystyle E_{n-1,i}&=f(e^{-\bar{c}_ihM}y_{n-1}+h\sum\limits_{j=1}^{i-1}\bar{a}_{ij}k_{n-1,j})
-f(e^{-\bar{c}_ihM}y(t_{n-1})+h\sum\limits_{j=1}^{i-1}\bar{a}_{ij}\bar{k}_{n-1,j}),\cr\noalign{\vskip1truemm}
\displaystyle e_{n}&=e^{-hM}e_{n-1}+h\sum\limits_{i=1}^{4}\bar{b}_i\big(k_{n-1,i}-\bar{k}_{n-1,i}\big)
+\bar{w}_4(z)-\hat{w}_4(z)-\delta_{n}.
\end{aligned}
\end{equation}
It follows from the first formula of (\ref{err1}) that
\begin{equation}\label{internal stages error}
\begin{aligned}
\displaystyle E_{l,1}&=f(y_l)-f(y(t_l)),\cr\noalign{\vskip2truemm}
\displaystyle E_{l,2}&=f(e^{-\bar{c}_2hM}y_{l}+h\bar{a}_{21}k_{l,1})
-f(e^{-\bar{c}_2hM}y(t_l)+h\bar{a}_{21}\bar{k}_{l,1}),\cr\noalign{\vskip2truemm}
\displaystyle E_{l,3}&=f(e^{-\bar{c}_3hM}y_{l}+h\bar{a}_{31}k_{l,1}+h\bar{a}_{32}k_{l,2})
-f(e^{-\bar{c}_3hM}y(t_l)+h\bar{a}_{31}\bar{k}_{l,1}+h\bar{a}_{32}\bar{k}_{l,2}),\cr\noalign{\vskip2truemm}
\displaystyle E_{l,4}&=f(e^{-\bar{c}_4hM}y_{l}+h\bar{a}_{41}k_{l,1}+h\bar{a}_{42}k_{l,2}
+h\bar{a}_{43}k_{l,3})-f(e^{-\bar{c}_4hM}y(t_l)+h\bar{a}_{41}\bar{k}_{l,1}+h\bar{a}_{42}\bar{k}_{l,2}\cr\noalign{\vskip2truemm}
&\quad +h\bar{a}_{43}\bar{k}_{l,3}).
\end{aligned}
\end{equation}
The  formula (\ref{internal stages error}) holds for $l=0,\ldots,n-1$. Inserting the Taylor series of $\bar{k}_{l,1},\ \bar{k}_{l,2},\ \bar{k}_{l,3},\ \bar{k}_{l,4}$ about $y(t_{l})$ into \eqref{residual},  then $\delta_{l+1}$ satisfies
\begin{equation*}
\begin{aligned}
\delta_{l+1}&=y(t_l)+hy'(t_l)+\frac{h^2}{2!}y''(t_l)+\frac{h^3}{3!}y'''(t_l)
+\frac{h^4}{4!}y^{(4)}(t_l)+\mathcal{O}(h^5)-e^{-hM}y(t_l)-h\big(\bar{b}_1\bar{k}_{l,1}\cr\noalign{\vskip2truemm}
&\quad+\bar{b}_2\bar{k}_{l,2}+\bar{b}_3\bar{k}_{l,3}+\bar{b}_4\bar{k}_{l,4}\big)-\hat{w}_4(z)\cr\noalign{\vskip1truemm}
=&y(t_l)+hg(t_l)+\frac{h^2}{2!}(-M+f'_y(y_l))g(t_l)+\frac{h^3}{3!}\Big(M^2g(t_l)
+(-M+f'_y(y(t_l)))f'_y(y(t_l))g(t_l)\cr\noalign{\vskip2truemm}
&-f'_y(y(t_l))Mg(t_l)+f''_{yy}(y(t_l))(g(t_l),g(t_l))\Big)
+\frac{h^4}{4!}\Big(-M^3g(t_l)+M^2f'_y(y(t_l))g(t_l)\cr\noalign{\vskip2truemm}
&-M f'_y(y(t_l))(-M+f'_y(y(t_l)))g(t_l)-Mf''_{yy}(y(t_l))(g(t_l),g(t_l))+f'''_{yyy}(y(t_l))(g(t_l),g(t_l),
\cr\noalign{\vskip2truemm}
&g(t_l))+3f''_{yy}(y(t_l))((-M+f'_y(y(t_l)))g(t_l),g(t_l))+f'_y(y(t_l))(-M+f'_y(y(t_l)))(-M\cr\noalign{\vskip2truemm}
&+f'_y(y(t_l)))g(t_l)+f'_y(y(t_l))f''_{yy}(y(t_l))(g(t_l),g(t_l))\Big) +\mathcal{O}(h^5)
 \cr\noalign{\vskip1truemm}
 & -\Bigg\{y(t_l)-hMy(t_l)+h(\bar{b}_1+\bar{b}_2+\bar{b}_3+\bar{b}_4)f(y(t_l))-\frac{h^2Mg(y(t_l))}{2!}+h^2(\bar{b}_2\bar{c}_{2}+\bar{b}_3\bar{c}_{3}\cr\noalign{\vskip2truemm}
&+\bar{b}_4\bar{c}_{4})f'_y(y(t_l))g(t_l)
-\frac{h^3}{3!}f'_y(y(t_l))Mf(y(t_l)) +\frac{h^3}{3!}M(M-f'_y(y(t_l))g(t_l)+\frac{h^3}{2!}
(\bar{b}_2\bar{c}_{2}^2\cr\noalign{\vskip2truemm}
 \end{aligned}
\end{equation*}
\begin{equation*}
\begin{aligned}
&+\bar{b}_3\bar{c}_{3}^2+\bar{b}_4\bar{c}_{4}^2)f''_{yy}(y(t_l))(g(t_l),g(t_l))+h^3(\bar{b}_3\bar{a}_{32}\bar{c}_{2}+\bar{b}_4(\bar{a}_{42}\bar{c}_{2}+\bar{a}_{43}\bar{c}_{3}))f'_y(y(t_l))\big(M^2y(t_l)\cr\noalign{\vskip4truemm}
&+f'_y(y(t_l))g(t_l)\big)+\frac{h^4}{4!} \Big(-M^3g(t_l)+M^2f'_y(y(t_l))g(t_l)
-Mf''_{yy}(y(t_l))(g(t_l),g(t_l))-M \cr\noalign{\vskip2truemm}
&\cdot f'_y(y(t_l))(-M+f'_y(y(t_l)))g(t_l)
-f'_y(y_l)Mf'_y(y(t_l))g(t_l)-f'_y(y(t_l))f'_y(y(t_l))
Mf(y(t_l))\cr\noalign{\vskip2truemm}
&+3f''_{yy}(y(t_l))(-Mf(y(t_l)),g(t_l))\Big)+h^4\bar{b}_4\bar{a}_{43}\bar{a}_{32}\bar{a}_{21}
f'_y(y(t_l))\big((M^2+f'_y(y(t_l))^2)g(t_l)\cr\noalign{\vskip2truemm}
&+f'_y(y(t_l))^2M^2y(t_l)\big)+\frac{h^4}{3!}(\bar{b}_2\bar{c}_{2}^3+\bar{b}_3\bar{c}_{3}^3
+\bar{b}_4\bar{c}_{4}^3) f'''_{yyy}(y(t_l))\big(g(t_l), g(t_l),g(t_l)\big)
\cr\noalign{\vskip2truemm}
&+h^4\big(\bar{b}_3\bar{c}_3\bar{a}_{32}\bar{a}_{21}+\bar{b}_4\bar{c}_4\bar{a}_{42}\bar{c}_2
+\bar{b}_4\bar{c}_4\bar{a}_{43}\bar{c}_3\big)f''_{yy}(y(t_l))\big(M^2y(t_l)+f'_y(y(t_l))g(t_l),g(t_l)\big)\cr\noalign{\vskip2truemm}
&+\frac{h^4}{2!}\big(\bar{b}_3\bar{a}_{32}\bar{c}_{2}^2+\bar{b}_4\bar{a}_{42}\bar{c}_{2}^2
+\bar{b}_4\bar{a}_{43}\bar{c}_{3}^2\big)f'_y(y(t_l))f''_{yy}(y(t_l))\big(g(t_l),g(t_l)\big)\Bigg\},
\end{aligned}
\end{equation*}
where $g(t_l)=-My(t_l)+f(y(t_l))$. Since the coefficients of the  explicit SVERK method satisfies the order conditions (\ref{order condition}),  one has
\begin{equation}
\delta_{l+1}=C_1h^5,
\end{equation}
where $C_1$ depends on $t_{end}$, $ \Vert M \Vert$, $f^{(p)}(t_n)$, $p=0,\ldots,4$, but is independent of
$n$ and $h$.

 Taking the norm $\Vert \cdot \Vert$ for the second formula of (\ref{err1}) and  (\ref{internal stages error})  yields
\begin{equation}\label{error}
\begin{aligned}
\Vert E_{l,1} \Vert &\leq L\Vert e_l \Vert,\cr\noalign{\vskip2truemm}
\Vert E_{l,2} \Vert &\leq L\Vert e^{-\bar{c}_2hM}e_l+h\bar{a}_{21}(k_{l,1}-\bar{k}_{l,1})\Vert
\leq L(1+h|\bar{a}_{21}| L)\Vert e_l \Vert,\cr\noalign{\vskip2truemm}
\Vert E_{l,3} \Vert &\leq L\Vert e^{-\bar{c}_3hM}e_l+h\bar{a}_{31}(k_{l,1}-\bar{k}_{l,1})
+h\bar{a}_{32}(k_{l,2}-\bar{k}_{l,2})\Vert\cr\noalign{\vskip2truemm}
&\leq L \big(1+h| \bar{a}_{31}| L+h|\bar{a}_{32}| L(1+h| \bar{a}_{21}| L)\big)\Vert e_l \Vert\cr\noalign{\vskip2truemm}
\Vert E_{l,4} \Vert 
&\leq L\Vert e^{-\bar{c}_4hM}e_l\Vert+Lh| \bar{a}_{41} | \cdot \Vert E_{l,1}  \Vert+
Lh| \bar{a}_{42} | \cdot  \Vert E_{l,2} \Vert + Lh | \bar{a}_{43} | \cdot  \Vert E_{l,3}\Vert \cr\noalign{\vskip2truemm}
&\leq L\big(1+h|\bar{a}_{41}|L+h|\bar{a}_{42}|L(1+h|\bar{a}_{21}|L)+h|\bar{a}_{43}|L(1+h|\bar{a}_{31} | L+h | \bar{a}_{32} | L\cr\noalign{\vskip2truemm}
&\quad +h^2| \bar{a}_{32}|  | a_{21} | L^2)\big) \Vert e_l \Vert
\end{aligned}
\end{equation}
and
$$\Vert e_{n} \Vert \leq  \Vert e_{n-1} \Vert +h\sum\limits_{i=1}^{4} | \bar{b}_i|
\cdot \Vert E_{n-1,i} \Vert + \Vert \bar{w}_4(z)-\hat{w}_4(z) \Vert+C_1h^5,$$

According to the expressions  of $\bar{w}_4(z)$ and $\hat{w}_{4}(z)$, one gets
\begin{equation}\label{residual2}
\begin{aligned}
\Vert \bar{w}_4(z)-\hat{w}_4(z) \Vert \leq C_2h^2L\Vert e_n \Vert +C_2h^3L\Vert e_n \Vert+C_2h^4L \Vert e_n \Vert,
\end{aligned}
\end{equation}
where $C_2$ depends on $ \Vert M \Vert$, $f^{(p)}(t_n)$, $p=0,1,2$.
It follows from  (\ref{error}) that
\begin{equation}\label{internal upper bound}
\begin{aligned}
h \sum\limits_{i=1}^4 |\bar{b}_i| \cdot \Vert  E_{l,i} \Vert &\leq  C_3hL\Vert e_l \Vert+C_3h^2L^2\Vert e_l \Vert+C_3h^3L^3\Vert e_l \Vert + C_3h^4L^4 \Vert e_l \Vert,
\end{aligned}
\end{equation}
for $l=0,\ldots,n-1,$ and the constant $C_3$ depends only on $a_{i,j}$ and $b_i$. We then obtain
\begin{equation}\label{updates}
\begin{aligned}
\Vert e_{n} \Vert &\leq \Vert e_{n-1} \Vert + h\sum\limits_{i=1}^4  |\bar{b}_i| \cdot \Vert E_{n-1,i} \Vert
+C_2h^2L \Vert e_{n-1} \Vert (1+h+h^2) +C_1h^5\\
&\leq \Vert e_1 \Vert +\sum\limits_{l=1}^{n-1} \big(h\sum\limits_{i=1}^4  |\bar{b}_i| \cdot \Vert E_{l,i} \Vert
+C_2h^2L \Vert e_{l} \Vert (1+h+h^2)\big)+(n-1)C_1h^5.
\end{aligned}
\end{equation}
Inserting (\ref{internal upper bound}) into (\ref{updates}) leads to
\begin{equation}\label{global error}
\Vert e_{n} \Vert \leq \sum\limits_{l=1}^{n-1} C_2C_3hL\big( 1+ hL 
+ h^2L^2 +h^3L^3 \big)\Vert e_l \Vert +C_1h^4.
\end{equation}
The application of the discrete Gronwall lemma to (\ref{global error}) gives
\begin{equation}
\Vert e_{n} \Vert \leq Ch^4,
\end{equation}
which completes the proof.
\end{proof}

The convergence of fourth-order explicit
MVERK methods can be analyzed in the same way which be skipped for brevity.

\section{Numerical Experiments}\label{sec5}
  In this section, we conduct some numerical experiments to illustrate  the accuracy and efficiency of  our ERK methods. Throughout the numerical experiments, the matrix-valued  $\varphi$-function of exponential integrators are evaluated by a  Pade approximant (see, e.g., \cite{Moler2003,Berland2007}), and the standard fourth-order ERK  method \eqref{Tableau-ERK51} with five stages  with a small time step  is employed as the reference solution.
We select the following fourth-order ERK methods for comparison:
 \begin{itemize}
  \item ERK41: the fourth-order  explicit ERK method \eqref{Tableau-ERK51} with five stages in \cite{Hochbruck2005a};
  \item ERK42: the fourth-order  explicit ERK method  \eqref{Tableau-ERK41} with four stages in \cite{Krogstad2005};
  \item MVERK41: the fourth-order  explicit MVERK method (\ref{MVERK41}) with four stages presented in this paper;
  \item MVERK42: the fourth-order   explicit MVERK method (\ref{MVERK42}) with four stages presented in this paper;
  \item SVERK41: the fourth-order  explicit SVERK method (\ref{SVERK41}) with four stages presented in this paper;
  \item SVERK42: the fourth-order  explicit SVERK method (\ref{SVERK42}) with four stages presented in this paper.
\end{itemize}

\textbf{Problem 1.}
We first consider the following averaged system in wind-induced
 oscillation (\cite{Mclachlan1998})
\begin{equation*}
\begin{aligned}& \left(
                   \begin{array}{c}
                     x_1 \\
                      x_2 \\
                   \end{array}
                 \right)
'= \left(
    \begin{array}{cc}
      -\zeta& -\lambda\\
       \lambda & -\zeta \\
    \end{array}
  \right)\left(
                   \begin{array}{c}
                     x_1 \\
                      x_2 \\
                   \end{array}
                 \right)+
\left(                                                         \begin{array}{c}                                             x_1x_2\\
\frac{1}{2}(x_1^2-x_2^2)                                  \end{array}
 \right),
\end{aligned}\end{equation*}
 where $\zeta \geq
0$ is a damping factor and $\lambda$ is a detuning parameter (see, e.g., \cite{Guck1983}).  Set
$$\zeta=r\cos(\theta),\ \lambda=r\sin(\theta),\ r \geq 0, \ 0 \leq \theta \leq \pi/2.$$
This system can be written as
\begin{equation*}
\begin{aligned}& \left(
                   \begin{array}{c}
                     x_1 \\
                      x_2 \\
                   \end{array}
                 \right)
'= \left(
    \begin{array}{cc}
      -\cos(\theta)& -\sin(\theta)\\
       \sin(\theta)& -\cos(\theta) \\
    \end{array}
  \right)
\left(
\begin{array}{c}
rx_1-\frac{1}{2}\sin(\theta)(x_2^2-x_1^2)-\cos(\theta)x_1x_2\\
rx_2-\sin(\theta)x_1x_2+\frac{1}{2}\cos(\theta)(x_2^2-x_1^2)
     \end{array}
       \right).
\end{aligned}\end{equation*}

\begin{figure}[!htb]
\centering
\begin{tabular}[c]{cccc}%
  \subfigure[]{\includegraphics[width=5.7cm,height=6.2cm]{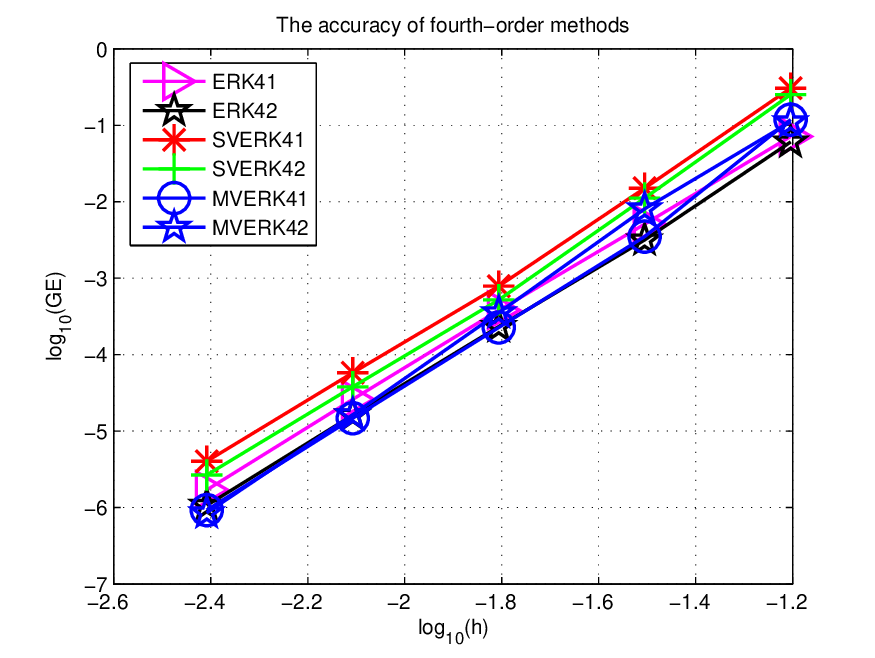}}
  \subfigure[]{\includegraphics[width=5.7cm,height=6.2cm]{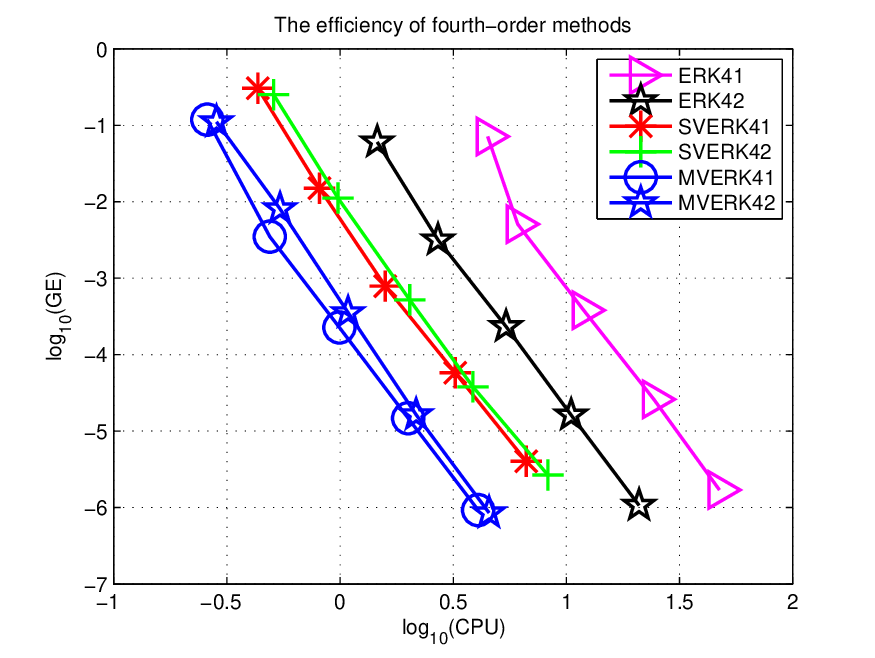}}
\end{tabular}
\caption{Results for  Problem 1. {(a)}: The $\log$-$\log$ plots of global
errors (GE) against $h$. {(b)}: The $\log$-$\log$ plots of global
errors against the CPU time.}\label{pro1}
\end{figure}

We integrate the conservative system  over the interval $[0,100]$  with  parameters $\theta=\pi/2$, $r=20$ and stepsizes $h=1/2^k$ for $k=4,\ldots,8$.
Fig. \ref{pro1} presents  the global errors $GE$
against the stepsizes and the CPU time  in a log–log scale. It can be observe that our methods have the
same accuracy with standard fourth-order  explicit exponential integrators, and the higher efficiency than standard exponential integrators  supported by their less CPU times (seconds).

\begin{figure}[!htb]
\centering
\begin{tabular}[c]{cccc}%
  \subfigure[]{\includegraphics[width=5.7cm,height=6.2cm]{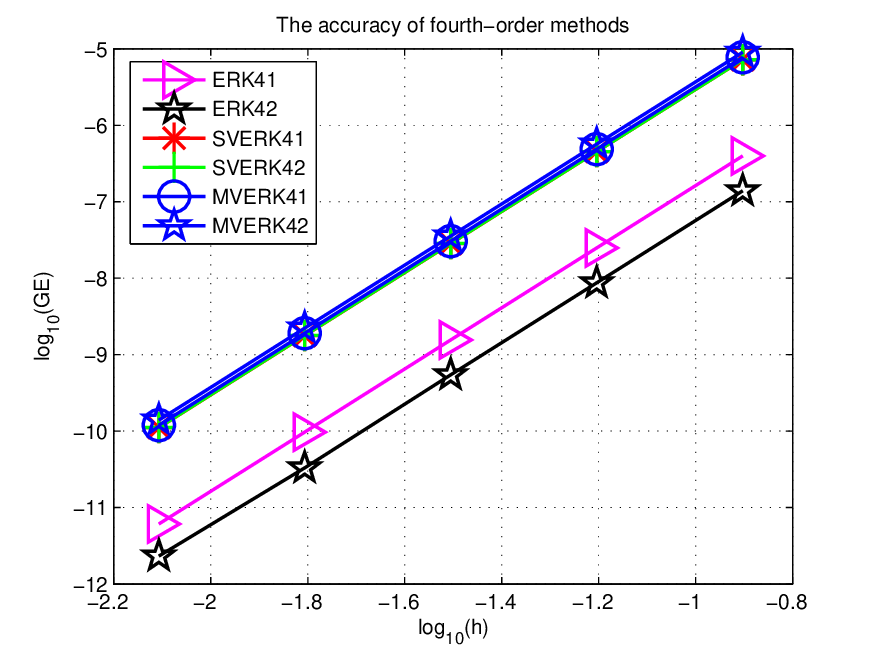}}
  \subfigure[]{\includegraphics[width=5.7cm,height=6.2cm]{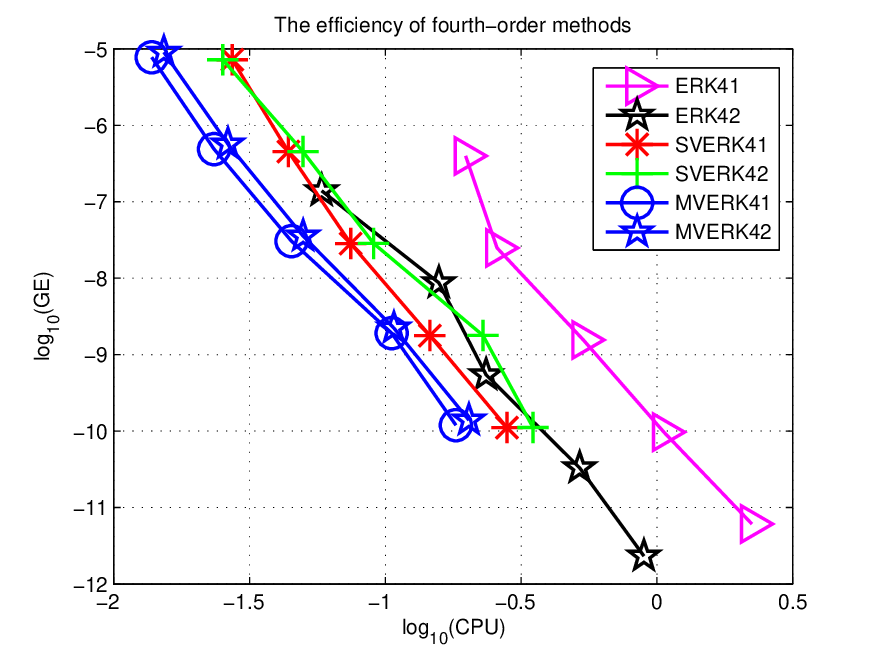}}
\end{tabular}
\caption{Results for  Problem 2. {(a)}: The $\log$-$\log$ plots of global
errors (GE) against $h$. {(b)}: The $\log$-$\log$ plots of global
errors against the CPU time.}\label{HHfirst}
\end{figure}

\textbf{Problem 2.}
 The H\'{e}non-Heiles Model  was used to describe the stellar motion (see, e.g., \cite{Hairer2006,Henon1964}), which has the following identical form
\begin{equation*}
\begin{aligned}& \left(
                   \begin{array}{c}
                     x_1 \\
                      x_2 \\
                      y_1 \\
                      y_2\\
                   \end{array}
                 \right)
'+\left(
    \begin{array}{cccc}
      0 & 0 & -1 &0\\
      0 & 0 & 0  &-1\\
      1 & 0 & 0  &0 \\
      0 & 1 & 0  &0\\
    \end{array}
  \right)\left(
                   \begin{array}{c}
                     x_1 \\
                      x_2 \\
                      y_1\\
                      y_2\\
                   \end{array}
                 \right)=
\left(
  \begin{array}{c}
  0\\
  0\\
  -2x_1x_2\\
  -x_1^2+x_2^2\\
  \end{array}
\right).
\end{aligned}
\end{equation*}
 We select the initial values as
 \begin{equation*}
\big(x_1(0),x_2(0),y_1(0),y_2(0)\big)^{\mathrm{T}}=(\sqrt{\frac{11}{96}},0,0,\frac{1}{4})^{\mathrm{T}}.
 \end{equation*}
 Fig. \ref{HHfirst} presents that  this problem is solved on the interval $[0,10]$  with stepsizes $h=1/2^{k},\ k=3,\ldots,7$ for ERK41, ERK42, SVERK41, SVERK42, MVERK41, and MVERK42.  

\begin{figure}[!htb]
\centering
\begin{tabular}[c]{cccc}%
  \subfigure[]{\includegraphics[width=5.7cm,height=6.2cm]{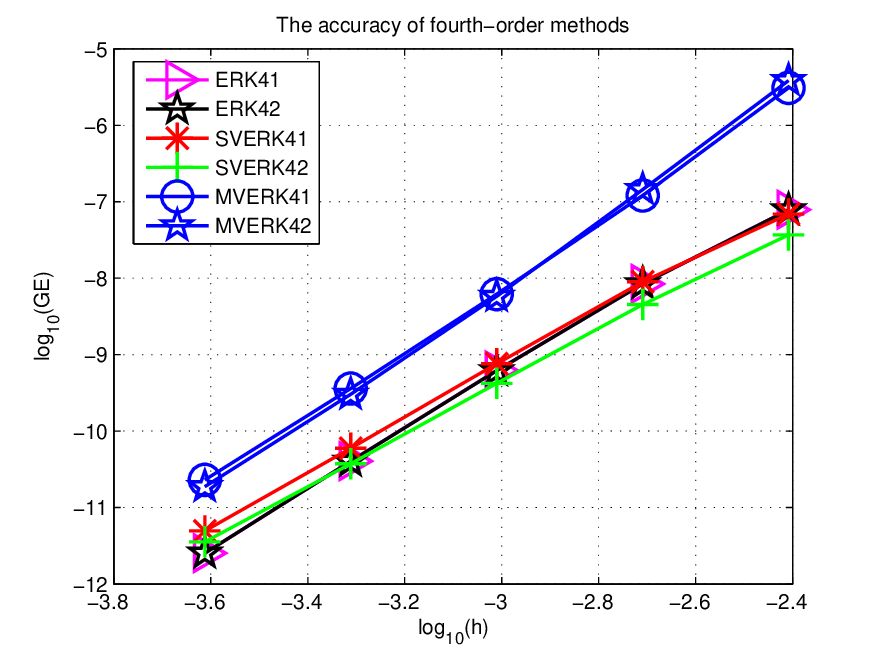}}
  \subfigure[]{\includegraphics[width=5.7cm,height=6.2cm]{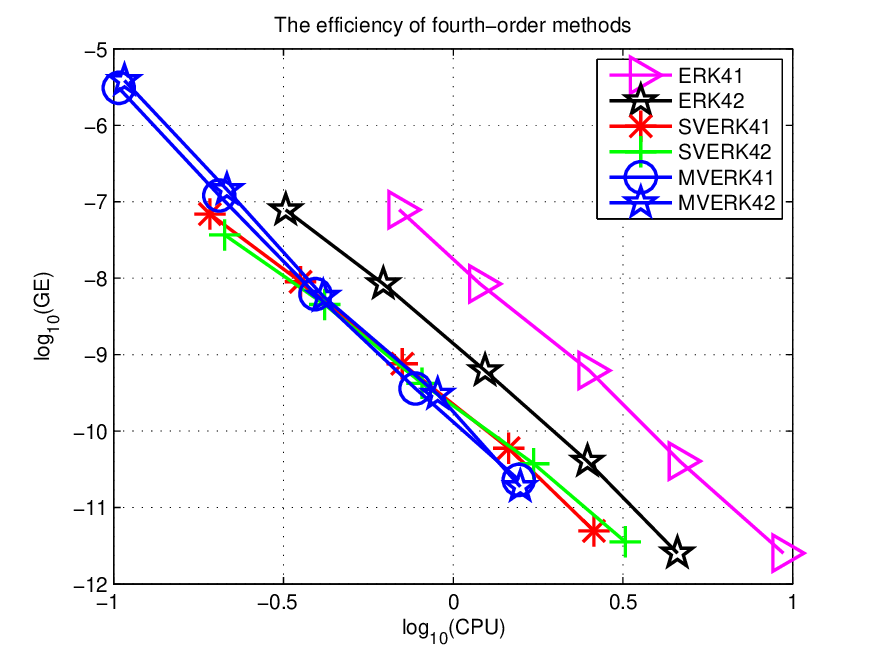}}
\end{tabular}
\caption{Results for  Problem 3. {(a)}: The $\log$-$\log$ plots of global
errors (GE) against $h$. {(b)}: The $\log$-$\log$ plots of global
errors against the CPU time.}\label{pro2}
\end{figure}

\textbf{Problem 3.}
We consider the Allen-Cahn equation \cite{Allen1979,feng}
\begin{equation*}
\left\{
\begin{aligned}
&\frac{\partial u}{\partial t}=\epsilon\frac{\partial^2u}{\partial x^2}+u-u^3, \quad -1<x<1,\ t>0,\\
&u(x,0)=0.53x+0.47\sin(-1.5\pi x),\quad -1\leq x \leq 1.
\end{aligned}\right.
\end{equation*}
  Allen and Cahn firstly introduced the equation  to describe the motion of anti-phase boundaries in crystalline solids \cite{Allen1979}.
After approximating the spatial derivatives with 32-point Chebyshev spectral method,
we obtain the following stiff system of first-order ordinary differential equations
\[\frac{dU}{dt}+MU=F(t,U), \quad t\in[0,t_{\rm end}],
\]
where $U(t)=(u_1(t),\cdots,u_N(t))^{\mathrm{T}},$ $u_i(t)\approx u(x_i,t),\
x_i=-1+i\Delta x$, for $i=1,{\ldots},N$ and $\Delta x=2/N$.
Here, the matrix $M$ is full and the nonlinear term is
$F(t,U)=u-u^3=(u_1-u_1^3,\cdots,u_N-u_N^3)^{\mathrm{T}}.$ We choose the
parameters $\epsilon=0.01,\  N=32$ and integrate the obtained stiff
system over the interval $[0,1]$. The numerical results are presented in Fig. \ref{pro2}
 with the stepsizes $h=1/2^k$ for $k=8,\ldots,12$.

\begin{figure}[!htb]
\centering
\begin{tabular}[c]{cccc}%
  \subfigure[]{\includegraphics[width=5.7cm,height=6.2cm]{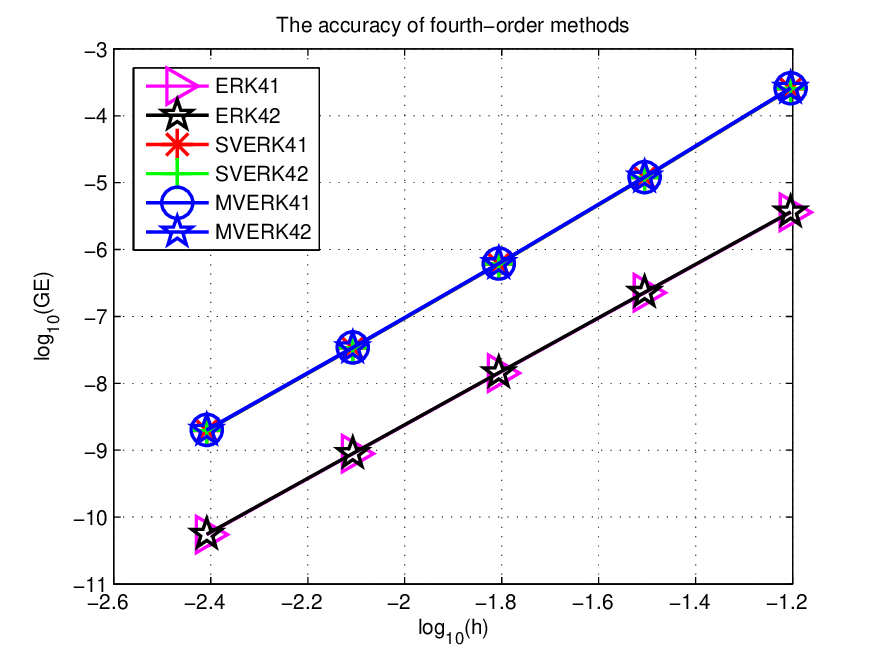}}
  \subfigure[]{\includegraphics[width=5.7cm,height=6.2cm]{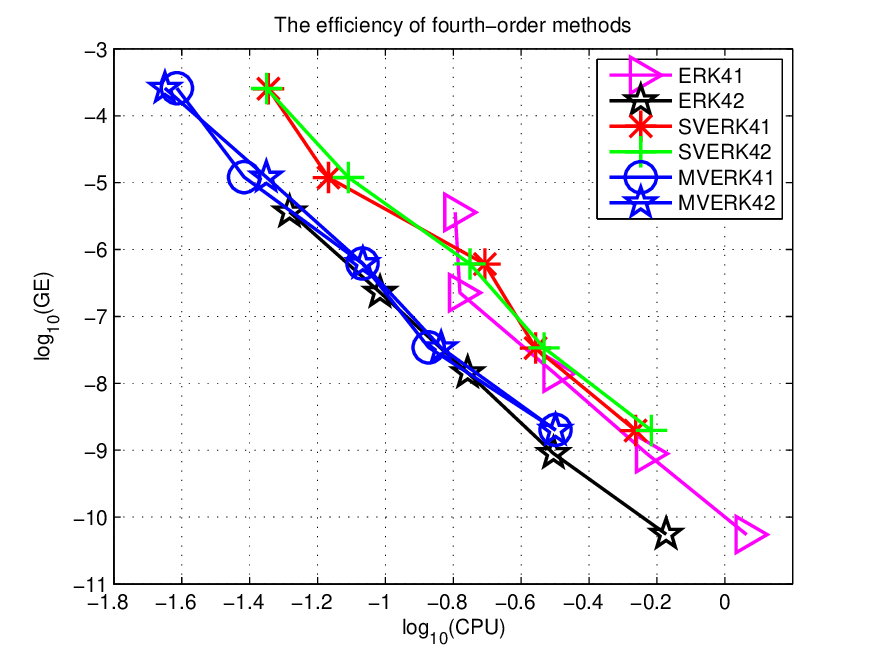}}
\end{tabular}
\caption{Results for  Problem 4. {(a)}: The $\log$-$\log$ plots of global
errors (GE) against $h$. {(b)}: The $\log$-$\log$ plots of global
errors against the CPU time.} \label{SG}
\end{figure}

\textbf{Problem 4.}\label{Sine-Gordon}  Consider the sine-Gorden equation with
periodic boundary conditions \cite{Fang2021,Franco2006}
\[\left\{\begin{array}{l}
\dfrac{\partial^2u}{\partial t^2}=\dfrac{\partial^2u}{\partial
x^2}-\sin(u), \quad -1<x<1, t>0,\cr\noalign{\vskip1truemm}
u(-1,t)=u(1,t).\cr\noalign{\vskip1truemm}
\end{array}\right.\]
Discretising the spatial derivative $\partial_{xx}$ by  the second-order symmetric differences yields
\[\frac{d}{dt}\left(\begin{array}[c]{c}
U^{\prime}\\
U\\
\end{array}\right)+\left(\begin{array}[c]{cc}
{\bf 0}&M\\
-I&{\bf 0}\\
\end{array}\right)\left(\begin{array}[c]{c}
U^{\prime}\\
U\\
\end{array}\right)=\left(\begin{array}[c]{c}
-\sin(U)\\
{\bf 0}\\
\end{array}\right), \quad t\in[0,t_{\rm end}].
\]
In here, $U(t)=(u_1(t),\cdots,u_N(t))^T$ with $u_i(t)\approx u(x_i,t)$
{for $i=1,\ldots,N$,}  with $\Delta x=2/N$ and $x_i=-1+i\Delta x$,
$F(t,U)=-\sin(u)=-(\sin(u_1),\cdots,\sin(u_N))^T$, and
\[M=\dfrac{1}{\Delta x^2}
\left(
\begin{array}
[c]{ccccc}
2 & -1 &  &  &-1 \\
-1 & 2&  -1 &  &  \\
& \ddots & \ddots & \ddots &   \\
  &  & -1 & 2 & -1\\
  -1&  &   &  -1&2  \\
\end{array}
\right).
\]
In this test, we choose the initial value conditions
\[U(0)=(\pi)^N_{i=1}, \ U'(0)=\sqrt{N}\left(0.01+\sin(\dfrac{2\pi i}{N})\right)_{i=1}^N\]
with $N=32$, and solve the problem on the interval $[0,1]$ with stepsizes $h=1/2^k,\ k=4,\ldots,8$.
The global errors GE against the stepsizes and the CPU time (seconds)
for ERK41, ERK42, SVERK41, SVERK42, MVERK41 and MVERK42  are respectively presented in  Fig. \ref{SG} (a) and (b).

\begin{figure}[!htb]
\centering
\begin{tabular}[c]{cccc}%
  \subfigure[]{\includegraphics[width=5.7cm,height=6.2cm]{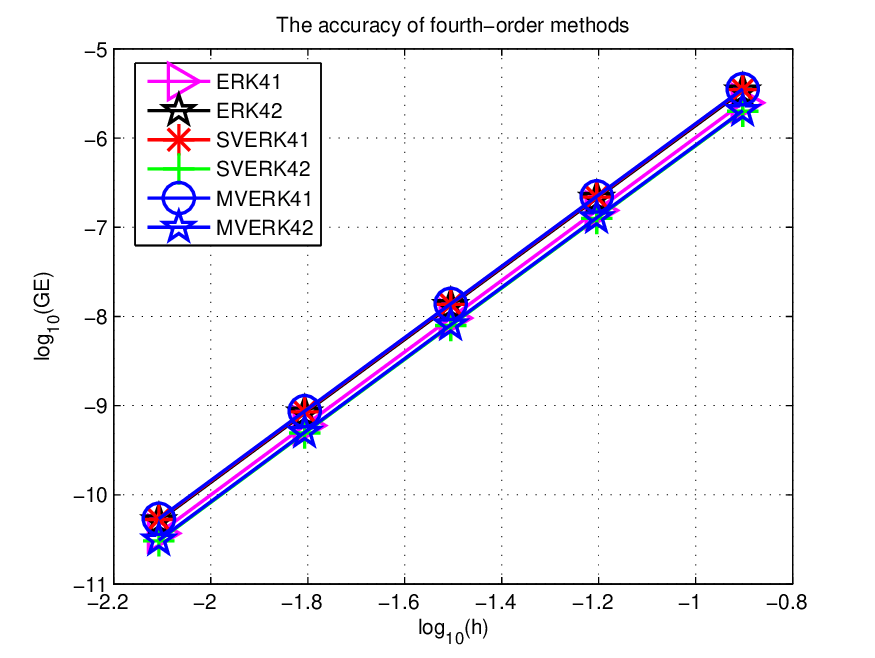}}
  \subfigure[]{\includegraphics[width=5.7cm,height=6.2cm]{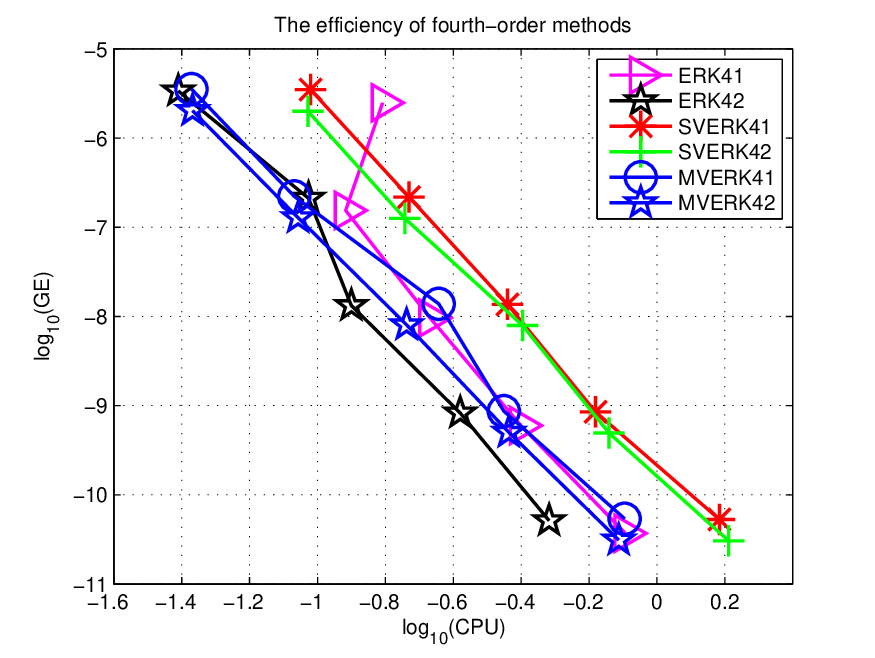}}
\end{tabular}
\caption{Results for  Problem 5. {(a)}: The $\log$-$\log$ plots of global
errors (GE) against $h$. {(b)}: The $\log$-$\log$ plots of global
errors against the CPU time.}\label{pro5}
\end{figure}

\textbf{Problem 5.}\label{Non-Schrodinger}
 We consider  the nonlinear
Schr\"{o}dinger equation (see  \cite{Chen2001})
\begin{equation*}\begin{aligned}
&i\psi_{t}+\psi_{xx}+2|\psi|^{2}\psi=0,\quad
 \psi(x,0)= 0.5 + 0.025 \cos(\mu x),\\
\end{aligned}
\end{equation*}
with the periodic boundary condition $\psi(0,t)=\psi(L,t).$
Letting $L =4\sqrt{2}\pi$ and $\mu =
2\pi/L$ and  $\psi = p + \textmd{i}q,$ we transform this equation into a
pair of real-valued equations
\begin{equation*}\label{semi}
\begin{aligned}
&p_t +q_{xx} + 2(p^2  + q^2)q = 0,\\
&q_t -p_{xx} -2(p^2  + q^2)p = 0.\\
\end{aligned}
\end{equation*}
Using the discretization on spatial variable with the
pseudospectral method leads to
\begin{equation}\label{semi sd}
\left(
  \begin{array}{c}
    \textbf{p} \\
    \textbf{q} \\
  \end{array}
\right)'=
 \begin{aligned}
 \left(
   \begin{array}{cc}
     0 & -D_2 \\
     D_2 & 0 \\
   \end{array}
 \right)\left(
  \begin{array}{c}
    \textbf{p} \\
    \textbf{q} \\
  \end{array}
\right)+\left(
          \begin{array}{c}
            -2(\textbf{p}^{2}+ \textbf{q}^{2})\cdot  \textbf{q} \\
            2( \textbf{p}^{2}+ \textbf{q}^{2})\cdot \textbf{p} \\
          \end{array}
        \right)
\end{aligned}
\end{equation}
where $\textbf{p}=(p_0,p_1,\ldots,p_{N-1})^{\mathrm{T}},\
\textbf{q}=(q_0,q_1,\ldots,q_{N-1})^{\mathrm{T}}$ and
$D_2=(D_2)_{0\leq j,k\leq N-1}$ is the pseudospectral differential
matrix defined by:
\begin{equation*}
(D_2)_{jk}=\left\{\begin{aligned}
&\frac{1}{2}\mu^2(-1)^{j+k+1}\frac{1}{\sin^2(\mu(x_j-x_k)/2)},\quad j\neq k,\\
&-\mu^2\frac{2(N/2)^2+1}{6},\quad\quad\quad\quad\quad\quad\ \  \   j=k.
\end{aligned}\right.
\end{equation*}
In this test, we choose $N=48$  and solve this problem
on $[0,1]$. The  global
errors $GE$  against the stepsizes and the CPU time are
stated in Fig. \ref{pro5} with the stepsizes $h=1/2^k$ for $k=4,\ldots,8$.

 All numerical results indicate that  MVERK methods  and SVERK methods have the comparable accuracy and efficiency in comparison with  standard 
exponential integrators, and the fourth-order explicit  MVERK methods is more efficient than SVERK methods due to the difference of $w_4(z)$ and $\bar{w}_4(z)$.

\section{Conclusions}\label{sec6}

 In this paper,  we present the fourth-order explicit MVERK methods and SVERK methods with four stages, and the order
conditions of these ERK methods are derived, which are exactly identical to the order conditions of explicit  RK methods. These ERK methods have the favorable property, which can exactly
integrate the linear system $y'(t)+My(t)=0$. Numerical results present the  comparable accuracy
and the lower computational cost of our fourth-order explicit ERK
methods.

\section*{Acknowledgements}

This research was supported by the National Natural Science Foundation of China (12071419).

\bibliographystyle{tfs}
\bibliography{ijnam}

\begin{thebibliography}{10}
\providecommand{\MR}{\relax\unskip\space MR }
\providecommand{\url}[1]{\normalfont{#1}}
\providecommand{\urlprefix}{Available at }

\bibitem{Abhulimen2014}
C.E. Abhulimen, \emph{Exponentially fitted third derivative three-step methods
  for numerical integration of stiff initial value problems}, Appl. Math.
  Comput. 243 (2014), pp. 446--453.
  \urlprefix\url{https://doi.org/10.1016/j.amc.2014.05.096}. \MR{3244492}

\bibitem{Allen1979}
S.M. Allen and J.W. Cahn, \emph{A microscopic theory for antiphase boundary
  motion and its application to antiphase domain coarsening}, Acta Metallurgica
  27 (1979), pp. 1085--1095.

\bibitem{Berland2005}
H.v. Berland, B. Owren, and B.r. Skaflestad, \emph{{$B$}-series and order
  conditions for exponential integrators}, SIAM J. Numer. Anal. 43 (2005), pp.
  1715--1727. \urlprefix\url{https://doi.org/10.1137/040612683}. \MR{2182146}

\bibitem{Berland2007}
H. Berland, B. Skaflestad, and W.M. Wright, \emph{Expint---a matlab package for
  exponential integrators}, ACM Trans. Math. Softw. 33 (2007), 4.

\bibitem{Butcher2008}
J.C. Butcher, \emph{Numerical methods for ordinary differential equations}, 2nd
  ed., John Wiley \& Sons, Ltd., Chichester, 2008,
  \urlprefix\url{https://doi.org/10.1002/9780470753767}. \MR{2401398}

\bibitem{Buvoli2022}
T. Buvoli and M.L. Minion, \emph{On the stability of exponential integrators
  for non-diffusive equations}, J. Comput. Appl. Math. 409 (2022), pp. Paper
  No. 114126, 17. \urlprefix\url{https://doi.org/10.1016/j.cam.2022.114126}.
  \MR{4383133}

\bibitem{Cash1981Siam}
J.R. Cash, \emph{Second derivative extended backward differentiation formulas
  for the numerical integration of stiff systems}, SIAM J. Numer. Anal. 18
  (1981), pp. 21--36. \urlprefix\url{https://doi.org/10.1137/0718003}.
  \MR{603428}

\bibitem{Cell2008}
E. Celledoni, D. Cohen, and B. Owren, \emph{Symmetric exponential integrators
  with an application to the cubic {S}chr\"{o}dinger equation}, Found. Comput.
  Math. 8 (2008), pp. 303--317.
  \urlprefix\url{https://doi.org/10.1007/s10208-007-9016-7}. \MR{2413146}

\bibitem{Chen2001}
J. Chen and M. Qin, \emph{Multi-symplectic {F}ourier pseudospectral method for
  the nonlinear {S}chr\"{o}dinger equation}, Electron. Trans. Numer. Anal. 12
  (2001), pp. 193--204. \MR{1847917}

\bibitem{Du2019}
Q. Du, L. Ju, X. Li, and Z. Qiao, \emph{Maximum principle preserving
  exponential time differencing schemes for the nonlocal {A}llen-{C}ahn
  equation}, SIAM J. Numer. Anal. 57 (2019), pp. 875--898.
  \urlprefix\url{https://doi.org/10.1137/18M118236X}. \MR{3945242}

\bibitem{Du2021}
Q. Du, L. Ju, X. Li, and Z. Qiao, \emph{Maximum bound principles for a class of
  semilinear parabolic equations and exponential time-differencing schemes},
  SIAM Rev. 63 (2021), pp. 317--359.
  \urlprefix\url{https://doi.org/10.1137/19M1243750}.

\bibitem{Eiermann2006}
M. Eiermann and O.G. Ernst, \emph{A restarted {K}rylov subspace method for the
  evaluation of matrix functions}, SIAM J. Numer. Anal. 44 (2006), pp.
  2481--2504. \urlprefix\url{https://doi.org/10.1137/050633846}. \MR{2272603}

\bibitem{Enright1974siam}
W.H. Enright, \emph{Second derivative multistep methods for stiff ordinary
  differential equations}, SIAM J. Numer. Anal. 11 (1974), pp. 321--331.
  \urlprefix\url{https://doi.org/10.1137/0711029}. \MR{351083}

\bibitem{Fang2021}
Y. Fang, X. Hu, and J. Li, \emph{Explicit pseudo two-step exponential
  {R}unge-{K}utta methods for the numerical integration of first-order
  differential equations}, Numer. Algorithms 86 (2021), pp. 1143--1163.
  \urlprefix\url{https://doi.org/10.1007/s11075-020-00927-4}. \MR{4211115}

\bibitem{feng}
X. Feng, H. Song, T. Tang, and J. Yang, \emph{Nonlinear stability of the
  implicit-explicit methods for the {A}llen-{C}ahn equation}, Inverse Probl.
  Imaging 7 (2013), pp. 679--695.
  \urlprefix\url{https://doi.org/10.3934/ipi.2013.7.679}. \MR{3105349}

\bibitem{Franco2006}
J.M. Franco, \emph{New methods for oscillatory systems based on {ARKN}
  methods}, Appl. Numer. Math. 56 (2006), pp. 1040--1053.
  \urlprefix\url{https://doi.org/10.1016/j.apnum.2005.09.005}. \MR{2234838}

\bibitem{Guck1983}
J. Guckenheimer and P. Holmes, \emph{Nonlinear oscillations, dynamical systems,
  and bifurcations of vector fields}, Applied Mathematical Sciences Vol.~42,
  Springer-Verlag, New York, 1983,
  \urlprefix\url{https://doi.org/10.1007/978-1-4612-1140-2}. \MR{709768}

\bibitem{Hairer2006}
E. Hairer, C. Lubich, and G. Wanner, \emph{Geometric numerical integration},
  2nd ed., Springer Series in Computational Mathematics Vol.~31,
  Springer-Verlag, Berlin, 2006, Structure-preserving algorithms for ordinary
  differential equations. \MR{2221614}

\bibitem{Henon1964}
M. H\'{e}non and C. Heiles, \emph{The applicability of the third integral of
  motion: {S}ome numerical experiments}, Astronom. J. 69 (1964), pp. 73--79.
  \urlprefix\url{https://doi.org/10.1086/109234}. \MR{158746}

\bibitem{Hochbruck1997}
M. Hochbruck and C. Lubich, \emph{On {K}rylov subspace approximations to the
  matrix exponential operator}, SIAM J. Numer. Anal. 34 (1997), pp. 1911--1925.
  \urlprefix\url{https://doi.org/10.1137/S0036142995280572}. \MR{1472203}

\bibitem{Hochbruck1998}
M. Hochbruck, C. Lubich, and H. Selhofer, \emph{Exponential integrators for
  large systems of differential equations}, SIAM J. Sci. Comput. 19 (1998), pp.
  1552--1574. \urlprefix\url{https://doi.org/10.1137/S1064827595295337}.
  \MR{1618808}

\bibitem{Hochbruck2005a}
M. Hochbruck and A. Ostermann, \emph{Explicit exponential {R}unge-{K}utta
  methods for semilinear parabolic problems}, SIAM J. Numer. Anal. 43 (2005),
  pp. 1069--1090. \urlprefix\url{https://doi.org/10.1137/040611434}.
  \MR{2177796}

\bibitem{Hochbruck2005b}
M. Hochbruck and A. Ostermann, \emph{Exponential {R}unge-{K}utta methods for
  parabolic problems}, Appl. Numer. Math. 53 (2005), pp. 323--339.
  \urlprefix\url{https://doi.org/10.1016/j.apnum.2004.08.005}. \MR{2128529}

\bibitem{Hochbruck2010}
M. Hochbruck and A. Ostermann, \emph{Exponential integrators}, Acta Numer. 19
  (2010), pp. 209--286.
  \urlprefix\url{https://doi.org/10.1017/S0962492910000048}. \MR{2652783}

\bibitem{Krogstad2005}
S. Krogstad, \emph{Generalized integrating factor methods for stiff {PDE}s}, J.
  Comput. Phys. 203 (2005), pp. 72--88.
  \urlprefix\url{https://doi.org/10.1016/j.jcp.2004.08.006}. \MR{2104391}

\bibitem{Lambert1972}
J.D. Lambert and S.T. Sigurdsson, \emph{Multistep methods with variable matrix
  coefficients}, SIAM J. Numer. Anal. 9 (1972), pp. 715--733.
  \urlprefix\url{https://doi.org/10.1137/0709060}. \MR{317548}

\bibitem{Lawson1967}
J.D. Lawson, \emph{Generalized {R}unge-{K}utta processes for stable systems
  with large {L}ipschitz constants}, SIAM J. Numer. Anal. 4 (1967), pp.
  372--380. \urlprefix\url{https://doi.org/10.1137/0704033}. \MR{221759}

\bibitem{Li2020}
B. Li, J. Yang, and Z. Zhou, \emph{Arbitrarily high-order exponential cut-off
  methods for preserving maximum principle of parabolic equations}, SIAM J.
  Sci. Comput. 42 (2020), pp. A3957--A3978.
  \urlprefix\url{https://doi.org/10.1137/20M1333456}. \MR{4186541}

\bibitem{Li2016}
Y. Li and X. Wu, \emph{Exponential integrators preserving first integrals or
  {L}yapunov functions for conservative or dissipative systems}, SIAM J. Sci.
  Comput. 38 (2016), pp. A1876--A1895.
  \urlprefix\url{https://doi.org/10.1137/15M1023257}.

\bibitem{Mclachlan1998}
R.I. McLachlan, G.R.W. Quispel, and N. Robidoux, \emph{Unified approach to
  {H}amiltonian systems, {P}oisson systems, gradient systems, and systems with
  {L}yapunov functions or first integrals}, Phys. Rev. Lett. 81 (1998), pp.
  2399--2403. \urlprefix\url{https://doi.org/10.1103/PhysRevLett.81.2399}.
  \MR{1643653}

\bibitem{Mei2022}
L. Mei, L. Huang, and X. Wu, \emph{Energy-preserving continuous-stage
  exponential {R}unge-{K}utta integrators for efficiently solving {H}amiltonian
  systems}, SIAM J. Sci. Comput. 44 (2022), pp. A1092--A1115.
  \urlprefix\url{https://doi.org/10.1137/21M1412475}. \MR{4417004}

\bibitem{Moler2003}
C. Moler and C. Van~Loan, \emph{Nineteen dubious ways to compute the
  exponential of a matrix, twenty-five years later}, SIAM Rev. 45 (2003), pp.
  3--49. \urlprefix\url{https://doi.org/10.1137/S00361445024180}. \MR{1981253}

\bibitem{Nosett1969}
S.P. N{\o}rsett, \emph{An {$A$}-stable modification of the {A}dams-{B}ashforth
  methods}, Lecture Notes in Math. Vol. Vol. 109, Springer, Berlin-New York,
  1969, pp. 214--219. \MR{267771}

\bibitem{Pope1963}
D.A. Pope, \emph{An exponential method of numerical integration of ordinary
  differential equations}, Comm. ACM 6 (1963), pp. 491--493.
  \urlprefix\url{https://doi.org/10.1145/366707.367592}. \MR{153118}

\bibitem{Qi2012}
R. Qi, C. Zhang, and Y. Zhang, \emph{Dissipativity of multistep {R}unge-{K}utta
  methods for nonlinear {V}olterra delay-integro-differential equations}, Acta
  Math. Appl. Sin. Engl. Ser. 28 (2012), pp. 225--236.
  \urlprefix\url{https://doi.org/10.1007/s10255-012-0142-x}. \MR{2914368}

\bibitem{Wang2022}
B. Wang, X. Hu, and X. Wu, \emph{Two new classes of exponential runge–kutta
  integrators for efficiently solving stiff systems or highly oscillatory
  problems}, International Journal of Computer Mathematics
  \urlprefix\url{https://doi.org/10.1080/00207160.2023.2294432}.

\bibitem{Wu2019}
B. Wang and X. Wu, \emph{Exponential collocation methods for conservative or
  dissipative systems}, J. Comput. Appl. Math. 360 (2019), pp. 99--116.
  \urlprefix\url{https://doi.org/10.1016/j.cam.2019.04.015}. \MR{3942760}

\bibitem{Wang2021}
W. Wang and C. Zhang, \emph{Dissipativity of variable-stepsize {R}unge-{K}utta
  methods for nonlinear functional differential equations with application to
  {N}icholson's blowflies models}, Commun. Nonlinear Sci. Numer. Simul. 97
  (2021), pp. Paper No. 105723, 19.
  \urlprefix\url{https://doi.org/10.1016/j.cnsns.2021.105723}. \MR{4209695}

\bibitem{Wu2018}
X. Wu and B. Wang, \emph{Recent developments in structure-preserving algorithms
  for oscillatory differential equations}, Science Press Beijing, Beijing;
  Springer, Singapore, 2018,
  \urlprefix\url{https://doi.org/10.1007/978-981-10-9004-2}. \MR{3791443}

\bibitem{Wu2021}
X. Wu and B. Wang, \emph{Geometric integrators for differential equations with
  highly oscillatory solutions}, Springer, Singapore; Science Press Beijing,
  Beijing, 2021, \urlprefix\url{https://doi.org/10.1007/978-981-16-0147-7}.
  \MR{4331421}

\bibitem{Zhang2014}
C. Zhang and T. Qin, \emph{The mixed {R}unge-{K}utta methods for a class of
  nonlinear functional-integro-differential equations}, Appl. Math. Comput. 237
  (2014), pp. 396--404.
  \urlprefix\url{https://doi.org/10.1016/j.amc.2014.03.143}. \MR{3201139}

\end{thebibliography}

\end{document}